\documentclass{amsart}

\usepackage{amssymb}
\usepackage{mathrsfs}
\usepackage{amsfonts}
\usepackage{verbatim}
\usepackage{pstricks}
\usepackage{pst-plot}

\newtheorem{thm}{Theorem}[section]
\newtheorem{lemma}[thm]{Lemma}
\newtheorem{defn}[thm]{Definition}
\newtheorem{rmk}[thm]{Remark}
\newtheorem{ex}[thm]{Example}
\newtheorem{cor}[thm]{Corollary}

\title{Geometric Structures on the Cochains of a Manifold}
\author{Scott O. Wilson}
\date{\today}

\begin{document}

\begin{abstract} 
In this paper we develop several algebraic structures on the simplicial 
cochains
of a triangulated manifold that are analogues of objects in differential 
geometry. 
We study a cochain product and prove several statements about its 
convergence to the wedge product on differential forms.
Also, for cochains with an inner product, we define a combinatorial 
Hodge star operator, and describe some applications, 
including a combinatorial period matrix for surfaces. 
We show that for a particularly nice cochain inner product,
these combinatorial structures converge to their continuum analogues as the
mesh of a triangulation tends to zero.
\end{abstract}

\maketitle


\section{Introduction} \label{sec:intro}
In this paper we develop combinatorial analogues of several objects 
in differential and complex geometry, including the Hodge star operator
and the period matrix of a Riemann surface. We define these 
structures on the appropriate combinatorial analogue of differential
forms, namely simplicial cochains.

As we recall in section~\ref{sec:smooth}, the two essential ingredients to the 
smooth Hodge star operator are Poincar\'e Duality and a metric, or inner 
product.
We'll define the combinatorial star operator in much the same way, using 
both an inner product and Poincar\'e Duality, expressed on cochains
in the form of a (graded) commutative product.

Using the inner product introduced in \cite{Do}, we prove the following:

\begin{thm} The combinatorial star operator, defined on the simplicial cochains
of a triangulated Riemannian manifold, converges to the smooth Hodge star
operator as the mesh of the triangulation tends to zero.
\end{thm}

We show in section~\ref{sec:surfaces} that, on a closed surface, 
this combinatorial star 
operator gives rise to a combinatorial period matrix and prove:

\begin{thm}
The combinatorial period matrix of a triangulated Riemannian 2-manifold 
converges to the conformal period matrix of the associated Riemann surface,
as the mesh of the triangulation tends to zero.
\end{thm}

This suggests a link between statistical mechanics and 
conformal field theory, where it is known that the partition function may 
be expressed in terms of theta functions of the conformal period 
matrix \cite{YM}, see also \cite{MS}, \cite{CM2}.

The above convergence statements are made precise by
using an embedding of simplicial cochains into differential forms, 
first introduced by Whitney \cite{Wh}. This approach was used
quite successfully by Dodziuk \cite{Do}, and later Dodziuk and Patodi 
\cite{DP}, to show that cochains provide a good approximation to 
smooth differential forms, and that the combinatorial Laplacian 
convergences to the smooth Laplacian. This formalism will be reviewed 
in section~\ref{sec:whitney}.

In section~\ref{sec:prod} we describe the cochain product that will be used in 
defining the combinatorial star operator. 
This product is of interest in its own right, and
we prove several results concerning its convergence to the wedge
product on forms; see also \cite{MK},\cite{CS}. 
These results may be of interest in
numerical analysis and the modeling of PDE's, since they give a computable
discrete model which approximates the algebra of smooth differential forms.
The convergence statements on the cochain product, theorems \ref{thm:prodconv} 
through ~\ref{thm:hmtps2zero}, are not needed for later sections.

In section~\ref{sec:star} we introduce the combinatorial star operator, 
and show that many of the interesting 
relations amongst $\star$, $d$, $\wedge$, and the adjoint  $d^{\ast}$ of $d$, 
that hold in the smooth setting, also hold in 
the combinatorial case. Some of the relations, though, are more illusive
and may only be recovered in the limit of a fine triangulation. 

In section~\ref{sec:surfaces} we study the combinatorial star operator 
on surfaces, and prove several results on the combinatorial period matrix, 
as mentioned above.

In the last two sections, ~\ref{sec:ip} and ~\ref{sec:circle}, 
we show how an explicit computation of the 
combinatorial star operator is related to ``summing over weighted
paths'', and perform these calculations for the circle.

I would like to thank Dennis Sullivan for his generosity, insight 
and many useful comments on this paper. 
I'd also like to thank Jozef Dodziuk for his help with several
points in his papers \cite{Do}, \cite{DP}, and
Ruben Costa-Santos and Barry McCoy for their inspiring work in \cite{MS}.

\section{Background and Acknowledgments} \label{sec:back}

In this section we describe previous results that are related to the contents 
of this paper. My sincere apologies to anyone whose work I have left out.

The cochain product we discuss in this paper was introduced by Whitney 
in \cite{Wh}.
It was also studied by Sullivan in the context of rational homotopy theory 
\cite{DPS}, by DuPont in his study of curvature and characteristic class 
\cite{JD}, and by Birmingham and Rakowski as a star product in 
lattice gauge theory \cite{BR}.
 
In connection with our result on the convergence of this cochain product to
the wedge product of forms, Kervaire has a related result for the 
Alexander-Whitney product $\cup$ on cochains \cite{MK}.
Kervaire states that, for differential forms $A, B$, and the associated 
cochains $a,b$, 
$$
\lim\limits_{k \to \infty} a \cup b \ (S^kc) = \int_c A \wedge B
$$
for a convenient choice of subdivisions $S^kc$ of the chain $c$.
Cheeger and Simons use this result in the context of cubical cell 
structures in \cite{CS}. 
There they construct an explicit map $E(A,B)$ satisfying
$$
\int A \wedge B - a \cup b = \delta E(A,B)
$$
and use it in the development of the theory of differential characters.
To the best of our knowledge, our convergence theorems for the 
commutative cochain product in section~\ref{sec:prod} are the first
to appear in the literature.

Several definitions of a discrete analogue of the Hodge-star operator
have been made. In \cite{MS}, Costa-Santos and McCoy define a 
discrete star operator for a particular 2-dimensional lattice and
study convergence properties as it relates to the Ising Model. Mercat
defines a discrete star operator for surfaces in \cite{CM}, using
a triangulation and its dual, and uses it to study a notion of discrete 
holomorphy and its relation to Ising criticality.

In \cite{TKB}, Tarhasaari, Kettunen and Bossavit describe how to make
explicit computations in electromagnetism using Whitney forms and 
a star operator defined using the de Rham map from forms to cochains. 
Teixeira and Chew have also defined Hodge operators on a lattices for
the purpose of studying electromagnetic theory.

Adams \cite{DA}, and also Sen, Sen, Sexton and Adams \cite{SSSA}, 
define two discrete star operators using a triangulation and its dual, and 
present applications to lattice gauge fields and Chern-Simons theory. 
de Beauc\'e and Sen \cite{BS} define star operators in a similar way and 
study applications to chiral Dirac fermions, and de Beauc\'e and Sen have 
generalized this to give a discreti\emph{z}ation scheme for differential 
geometry \cite{BS2}.

In the approaches using a triangulation and its dual, the star operator(s) 
are formulated using the duality map between the two cell decompositions. 
This map yields Poincar\'e Duality on (co)homology.
In this paper, we express Poincar\'e Duality by a commutative cup 
product on cochains and combine it with a non-degenerate inner product 
to define the star operator. Working this way, we obtain a single operator 
from one complex to itself.

Our convergence statements in section~\ref{sec:star}
are proven using the inner product 
introduced in \cite{Do}, and to the best of our knowledge, these are the
first results proving a convergence theorem for a discrete analogue
of the Hodge star operator.

In Dodziuk's paper \cite{Do}, and in \cite{DP} by 
Dodziuk and Patodi, the authors study a combinatorial Laplacian
on the cochains and proved that its eigenvalues converge to the smooth 
Laplacian. Such discrete notions of a Hodge structure, 
along with finite element method techniques, were used by Kotiuga 
\cite{RK}, and recently by Gross and Kotiuga \cite{GK}, in
the study of computational electromagnetism. Jin has used related techniques 
in studying electrodynamics \cite{JJ}. 

In connection with our application of the combinatorial star operator to 
surfaces, in particular proving the convergence of our
combinatorial period matrix to the conformal period matrix, Mercat
has a related result in \cite{CM2}. As part of his extensive study of 
what he calls ``discrete Riemann surfaces'', he assigns to any such 
object a ``period matrix'' of twice the expected size. 
He shows that there are two sub-matrices of the appropriate dimension 
($g \times g$) satisfying the property that, 
given what he calls ``a refining sequence of critical maps'', they both 
converge to the continuum period matrix of an associated 
Riemann surface. This uses his results on discrete holomorphy approximations
presented in \cite{CM3}. As with the star operator above, our approach
differs in that there is no ``doubling'' of complexes or operators.

There is another discussion of discrete period matrices presented in 
\cite{GY}. There Xianfeng Gu and Shing-Tung Yau give explicit 
algorithms for computing a period matrix for a surface. 
They point out that these can be implemented on the simplicial 
cochains by the use of the integration map from piecewise linear forms to 
simplicial cochains.

\section{Smooth Setting} \label{sec:smooth}
We begin with a brief review of some elementary definitions. Let $M$ be 
a closed oriented Riemannian $n$-manifold. 
A Riemannian metric induces an inner product on 
$\Omega(M)= \bigoplus_{j} \Omega^{j} = 
\bigoplus_{j} \Gamma(\bigwedge^{j}T^{*}M)$ in the following way: 
a Riemannian metric determines an inner product on 
$T^{\ast}M_{p}$ for all $p$, and hence an inner product for each $j$
on $\bigwedge^{j} T^{\ast}M_{p}$ (explicitly, via an orthonormal basis).
An inner product $\langle , \rangle$ on $\Omega (M)$ is then obtained by 
integration over M. 
If we denote the induced norm on $\bigwedge^{j} T^{\ast}M_{p}$ by 
$\lvert \hspace{1em} \rvert_{p}$, then the norm $\lVert \hspace{1em} \rVert$ 
on $\Omega(M)$ is given by
\begin{displaymath}
\lVert \omega \rVert 
= \Big( \int_{M} \lvert \omega \rvert_{p}^{2} \ dV\Big)^{1/2}
\end{displaymath}
where $dV$ is the Riemannian volume form on $M$.

Let $\mathcal{L}_{2} \Omega(M)$ denote the completion of $\Omega(M)$ 
with respect to this norm. We also use $\lVert \hspace{1em} \rVert$ to 
denote the norm on the completion.
Let the exterior derivative $d:\Omega^{j}(M) \to \Omega^{j+1}(M)$ be 
defined as usual.

\begin{defn} The Poincare-Duality pairing $(,):\Omega^{j}(M) \otimes 
\Omega^{n-j}(M) \to \mathbb{R}$ is defined by:
\begin{displaymath}
(\omega,\eta)= \int_{M} \omega \wedge \eta.
\end{displaymath}
\end{defn}

The pairing $(,)$ is bilinear, (graded) skew-symmetric and non-degenerate. 
It induces an isomorphism $(,):\Omega^{j}(M) \to (\Omega^{n-j}(M))^{\ast}$, 
where here $\ast$
denotes linear the dual. The induced map $\langle,\rangle:\Omega^{n-j}(M) \to 
(\Omega^{n-j}(M))^
{\ast}$ is also an isomorphism and one may check that the composition 
$\langle,\rangle^{-1} \circ ( , )$ equals the following operator:

\begin{defn} \label{defn:star}
The Hodge star operator $\star:\Omega^{j}(M) \to \Omega^{n-j}(M)$
is defined by:
\begin{displaymath}
\langle \star \omega,\eta \rangle = (\omega, \eta).
\end{displaymath}
\end{defn}

One may also define the operator $\star$ using local coordinates, see Spivak
\cite{Sp}. We note that this approach and the former definition give 
rise to the same operator $\star$
on $\mathcal{L}_{2} \Omega(M)$. We prefer to emphasize 
definition \ref{defn:star} since
it motivates definition ~\ref{defn:bigstar}, the combinatorial star operator.

\begin{defn} The adjoint of $d$, denoted by $d^{\ast}$, is defined by 
$\langle d^{\ast} \omega,\eta \rangle = \langle \omega, d \eta \rangle$. 
\end{defn}

Note that $d^{\ast}: \Omega^{j}(M) \to \Omega^{j-1}(M)$.
The following relations hold among $\star$, $d$ and $d^{\ast}$. See Spivak 
\cite{Sp}.

\begin{thm} As maps from $\Omega^{j}(M)$ to their respective ranges:
\begin{enumerate}
\item $\star d = (-1)^{j+1} d^{\ast} \star$
\item $\star d^{\ast} = (-1)^{j} d \star$
\item $\star^{2} = (-1)^{j(n-j)}Id$
\end{enumerate}
\end{thm}

\begin{defn} The Laplacian is defined to be $\Delta = d^{\ast} d + d d^{\ast}$.
\end{defn}

Finally, we state the Hodge decomposition theorem for $\Omega(M)$.
Let $\mathcal{H}^{j}(M) = \{ \omega \in \Omega^{j}(M) | \Delta \omega = 
d \omega = d^{\ast} \omega = 0 \}$ be the space of harmonic $j$-forms.

\begin{thm} There is an orthogonal direct sum decomposition
\begin{displaymath}
\Omega^{j}(M) \cong d \Omega^{j-1}(M) \oplus \mathcal{H}^{j}(M) 
\oplus d^{\ast} \Omega^{j+1}(M)
\end{displaymath}
and $\mathcal{H}^{j}(M) \cong H^{j}_{DR}(M)$, the De Rham cohomology of $M$
in degree $j$.
\end{thm} 

\section{Whitney Forms} \label{sec:whitney}

In his book, `Geometric Integration Theory', Whitney explores the idea 
of using cochains as integrands \cite{Wh}. A main result is
that such objects provide a reasonable integration theory that 
in some sense generalize the smooth theory of integration of 
differential forms.
This idea has been made even more precise by the work of Dodziuk \cite{Do},
who used a linear map of cochains into $\mathcal{L}_{2}$-forms
(due to Whitney \cite{Wh}) to show that cochains provide a good
approximation of differential forms. 
In this section we review some of the results.
The techniques involved illustrate a tight (and analytically precise) 
connection 
between cochains and forms, and will be used later to give precise meaning 
to our constructions on cochains. In particular, all of our convergence
statements about combinatorial and smooth objects will be cast 
in a similar way.

Let $M$ be a fixed closed smooth $n$-manifold and $K$ a fixed $C^{\infty}$ 
triangulation of $M$. We identify $K$ 
and $M$ and fix an ordering of the vertices of $K$. Let $C^{j}$ denote the 
simplicial cochains of degree $j$ 
of $K$ with values in $\mathbb{R}$. Given an ordering of the vertices of $K$,
we have a coboundary operator $\delta:C^{j} \to C^{j+1}$. 
Let $\mu_{i}$ denote the barycentric coordinate corresponding to the $i^{th}$ 
vertex $p_{i}$ of $K$. 
Since $M$ is compact, we may identify the cochains and chains of $K$ and 
for $c \in C^{j}$ write 
$c = \sum_{\tau} c_{\tau} \cdot \tau$ where $c_{\tau} \in \mathbb{R}$ and 
the sum over all $j$-simplicies 
$\tau = [p_{0},p_{1}, \dots ,p_{j}]$ of $K$ whose vertices 
form an increasing sequence 
with respect to the ordering of vertices in $K$.
We now define the Whitney embedding of cochains into 
$\mathcal{L}_{2}$-forms:

\begin{defn} For $\tau$ as above, we define
\[
W \tau = j! \sum_{i=0}^{j} (-1)^{i} \mu_{i} \ d \mu_{0} 
\wedge \dots \wedge 
\widehat{d \mu_{i}}
\wedge \dots \wedge d \mu_{j}.
\]
W is defined on all of $C^{j}$ by extending linearly.
\end{defn}

Note that the coordinates $\mu_{\alpha}$ are not even of class $C^{1}$, 
but they are $C^{\infty}$ on the interior
of any $n$-simplex of $K$. Hence, $d \mu_{\alpha}$ is defined 
and $W\tau$ is a well defined element of
$\mathcal{L}_{2} \Omega^{j}$. By the same consideration, $dW$ is also well 
defined.
Note both sides of the definition of $W$ are alternating, so this map is well 
defined for all simplicies regardless of the ordering of vertices.

Several properties of the map $W$ are given below. See \cite{Wh},\cite{Do},
\cite{DP} for details. 

\begin{thm} The following hold:
\begin{enumerate}
\item $W\tau = 0$ on $M \backslash \overline{St(\tau)}$
\item $dW=W\delta$
\end{enumerate}
where $St$ denotes the open star and $\overline{\ \ }$ denotes closure.
\end{thm}

One also has a map $R:\Omega^{j}(M) \to C^{j}(K)$, the de Rham map, given by
integration. Precisely, for any differential form $\omega$ and chain $c$ 
we have:
\begin{displaymath}
R\omega(c) = \int_{c} \omega
\end{displaymath}

It is a theorem of de Rham that this map is a quasi-isomorphism 
(it is a chain map 
by Stokes Theorem). $RW$ is well defined and one can check that $RW = Id$, 
see \cite{Wh}, \cite{Do}, \cite{DP}.

Before stating Dodziuk and Patodi's theorem that $WR$ is approximately equal 
to the identity, we first
give some definitions concerning triangulations. They also appear \cite{DP}.

\begin{defn} Let $K$ be a triangulation of an $n$-dimensional manifold $M$. 
The mesh $\eta = \eta(K)$ of a triangulation is: 
\begin{displaymath}
\eta = \sup r(p,q),
\end{displaymath}
where $r$ means the geodesic distance in $M$ and 
the supremum 
is taken over all pairs of vertices $p$, $q$ of a $1$-simplex in $K$.

The fullness $\Theta = \Theta(K)$  of a triangulation $K$ is  
\begin{displaymath}
\Theta(K) = \inf \frac{vol (\sigma)}{\eta^n},
\end{displaymath}
where the $\inf$ is taken over all $n$-simplexes $\sigma$ of $K$ and 
$vol (\sigma)$ is the Riemannian volume of $\sigma$, as a 
Riemannian submanifold of $M$.
\end{defn}

A Euclidean analogue of the following lemma was proven by Whitney in 
\cite{Wh} (IV.14). 
\begin{lemma} \label{lemma:techass} 
Let $M$ be a smooth Riemannian $n$-manifold.
\begin{enumerate}
\item Let $K$ be a smooth triangulation of $M$. Then there is a 
positive constant $\Theta_{0} > 0$ and a sequence of subdivisions
$K_1, K_2, \dots$ of K such that $\lim_{n \to \infty} \eta(K_n) = 0$ and
$\Theta (K_n) \geq \Theta_0$ for all n.
\item Let $\Theta_0 > 0$. There exist positive constants $C_1,C_2$ depending
on $M$ and $\Theta_0$ such that for all smooth triangulations $K$ of $M$ 
satisfying $\Theta(K) \geq \Theta_0$, all n-simplexes of 
$\sigma = [ p_0 , p_1 , \dots , p_n ]$ and vertices $p_k$ of $\sigma$,
{\setlength\arraycolsep{2pt}
\begin{eqnarray*}
vol(\sigma) &\leq& C_1 \cdot \eta^n \\
C_2 \cdot \eta &\leq& r( p_k , \sigma_{p_k}), 
\end{eqnarray*}}
where $r$ is the Riemannian distance, $vol(\sigma)$ is the Riemannian volume,
and $\sigma_{p_k} = [ p_0 ,\dots , p_{k-1} , p_{k+1} , \dots , p_n ]$ is the 
face of $\sigma$ opposite to $p_k$.
\end{enumerate}
\end{lemma}
Since any two metrics on $M$ are commensurable, the lemma follows from 
Whitney's Euclidean result, see also \cite{DP}.

We consider only those triangulations with fullness bounded 
below by some positive real constant $\Theta_{0}$. 
By the lemma, this guarantees that the volume of a simplex is on the order 
of it's mesh raised to the power of its dimension.
Geometrically, this means that in a sequence of triangulations, the 
shapes do not become too thin. (In fact, Whitney's \emph{standard subdivisions}
yield only finitely many shapes, and can be used to prove the 
first part of the lemma.) Most of the estimates in this paper depend on 
$\Theta_0$, as can be seen in the proofs. We'll not indicate this
dependence in the statements.

The following theorems are proved by Dodziuk and Patodi in \cite{DP}. They 
show that for a fine triangulation,
$WR$ is approximately equal to the identity. In this sense
the theorems give precise meaning to the statement: for a fine
triangulation, cochains provide a good approximation to differential 
forms.

\begin{thm} \label{thm:app} Let $\omega$ be a smooth form on $M$, and $\sigma$
be an n-simplex of $K$. There exists a constant $C$, independent of f, K and 
$\sigma$, such that
$$
\lvert \omega(p) - W R \omega \rvert_p  \leq C \cdot sup \bigg| 
\frac{\partial \omega}{\partial x^i} \bigg| \cdot \eta
$$
for all $p \in \sigma$. The supremum is taken over all $p \in \sigma$ and 
$i = 1, 2, \dots n$, and the partial derivatives are taken with
respect to a coordinate neighborhood containing $\sigma$.
\end{thm}

\begin{proof} Using similar techniques we'll prove a generalization; 
see Theorem ~\ref{thm:prodconv} and Remark ~\ref{rmk:prodrmk}.
\end{proof}

By integrating the above point-wise and applying a Sobolev inequality, Dodziuk
and Patodi \cite{DP} obtain the following:

\begin{cor} \label{cor:app} There exists a positive constant $C$ and a 
positive integer $m$, independent of $K$, such that
\[
\lVert \omega - WR \omega \rVert \le C \cdot 
\lVert (Id + \Delta)^{m} \omega \rVert \cdot \eta
\]
for all $C^{\infty}$ $j$-forms $\omega$ on $M$.
\end{cor}

\begin{proof} This is a special case of Corollary ~\ref{cor:prodconv2}.
\end{proof}

Now suppose the cochains $C$ are equipped with a non-degenerate
inner product $\langle , \rangle$ such
that, for distinct $i,j$, $C^i$ and $C^j$ are orthogonal. Then one can 
can define further structures on the cochains. In particular we have the 
following

\begin{defn} The adjoint of $\delta$, denoted by $\delta^{\ast}$, is 
defined by $\langle \delta^{\ast} \sigma, \tau \rangle =
\langle \sigma, \delta \tau \rangle$.
\end{defn}

Note that $\delta^{\ast}: C^{j}(K) \to C^{j-1}(K)$ is also squares to zero. 
One can also define

\begin{defn} The combinatorial Laplacian is defined to be 
$\blacktriangle = \delta^{\ast} \delta + \delta \delta^{\ast}$.
\end{defn}

Clearly, both $\delta^{\ast}$ and $\blacktriangle$ depend upon the 
choice of inner product. For any choice of non-degenerate inner
product, these operators give rise to 
a combinatorial Hodge theory: the space of harmonic $j$-cochains 
of $K$ is defined to be
$$\mathcal{H}C^{j}(K) = \{ a \in C^{j} | \blacktriangle a = \delta a 
= \delta^{\ast} a = 0 \}.
$$

The following theorem is due to Eckmann \cite{EC}:

\begin{thm} Let $C$ be a cochain complex with inner product 
$\langle , \rangle$, and induced differentials $\delta$ and $\delta^{\ast}$
as above. There is an orthogonal direct sum decomposition
\begin{displaymath}
C^{j}(K) \cong \delta C^{j-1}(K) \oplus \mathcal{H}C^{j}(K) \oplus 
\delta^{\ast} C^{j+1}(K)
\end{displaymath}
and $\mathcal{H}C^{j}(K) \cong H^{j}(K)$, the cohomology of $(K,\delta)$
in degree $j$.
\end{thm}

\begin{proof} We'll write $C^j$ for $C^j(K)$. The second statement of the
theorem follows from the first.

Using the fact that $\delta^{\ast}$ is the adjoint of $\delta$,
so $\delta \delta = \delta^{\ast} \delta^{\ast} = 0$, it is easy to check that 
$\delta C^{j-1}$, $\mathcal{H}C^{j}$, and $\delta^{\ast} C^{j+1}$
are orthogonal. Thus, it suffice to show:
$$
dim(C^{j}) = dim(\delta C^{j-1} \oplus \mathcal{H}C^{j} \oplus 
\delta^{\ast} C^{j+1})
$$
Let $\delta^{\ast}_j$ denote $\delta^{\ast}$ restricted to $C^j$. 
By orthogonality we have:
$$
dim(C^{j}) - dim(\delta^{\ast} C^j)
= dim(Ker(\delta^{\ast}_j)) 
= dim(\mathcal{H}C^{j}) + dim (\delta^{\ast} C^{j+1}).
$$

The proof is complete by showing 
$dim(\delta^{\ast} C^{j}) = dim(\delta C^{j-1})$.
This follows since, by the adjoint property, both 
$\delta: \delta^{\ast} C^{j} \to \delta C^{j-1}$ and 
$\delta^{\ast}: \delta C^{j-1} \to \delta^{\ast} C^{j}$
are injections of finite dimensional vector spaces.
\end{proof}

If $K$ is a triangulation of a Riemannian manifold $M$, then there 
is a particularly nice inner product on $C(K)$, which we'll call the 
\emph{Whitney inner product}. It is induced by the metric 
$\langle,\rangle$ on $\Omega(M)$ and the Whitney embedding of cochains
into $\mathcal{L}_2$-forms. We'll use the same notation $\langle,\rangle$ 
for this pairing on $C$: $\langle \sigma, \tau \rangle = 
\langle W \sigma, W \tau \rangle$. 

It is proven in \cite{Do} that the Whitney inner product on $C$ is 
non-degenerate. Further consideration of this inner product will be given 
in later sections. 
For now, following \cite{Do} and \cite{DP}, we describe how the combinatorial
Hodge theory, induced by the Whitney inner product, is related to the 
smooth Hodge theory. Precisely, we have the following theorem due to 
Dodziuk and Patodi \cite{DP}, which shows that the approximation 
$WR \approx Id$ respects the Hodge decompositions of $\Omega(M)$ and $C(K)$.

\begin{thm} \label{thm:happ}
Let $\omega \in \Omega^{j}(M)$, $R\omega \in C^{j}(K)$ have Hodge 
decompositions
{\setlength\arraycolsep{2pt}
\begin{eqnarray*}
\omega & = & d \omega_{1} + \omega_{2} + d^{\ast} \omega_{3}  \\ 
R\omega & = & \delta a_{1} + a_2 + \delta^{\ast} a_3
\end{eqnarray*}}
Then,
{\setlength\arraycolsep{2pt}
\begin{eqnarray*}
\| \omega_{1} - Wa_{1} \| & \leq & \lambda \cdot 
\| (Id + \Delta)^{m} \omega\| \cdots \eta \\
\| d\omega_{2} - W \delta a_{2} \| & \leq & \lambda \cdot
\| (Id + \Delta)^{m} \omega\| \cdot \eta \\
\| d^{\ast} \omega_{3} - W \delta^{\ast} a_{3} \| & \leq & \lambda 
\cdot \| (Id + \Delta)^{m} \omega \cdot \eta \| 
\end{eqnarray*}}
where $\lambda$ and $m$ are independent of  $\omega$ and $K$.
\end{thm}

\section{Cochain Product} \label{sec:prod}

In this section we describe a commutative, but non-associative, 
cochain product. It is of interest in its own right, and will be used 
to define the combinatorial star operator.

The product we define is induced by the Whitney embedding 
and the wedge product on forms, but also has
a nice combinatorial description. An easy way to state this is as follows: 
the product of a 
$j$-simplex and $k$-simplex is zero unless these simplicies are faces of a 
$(j+k)$-simplex, in which case the product 
is a rational multiple of this $(j+k)$-simplex.
We will prove a convergence theorem for this product, and also show that this
product's deviation from being associative converges to zero for 
`sufficiently smooth' cochains.

From the point of homotopy theory, it is natural to consider this 
commutative cochain product as part of a C-infinity algebra. 
We use Sullivan's local construction of a C-infinity algebra \cite{DS},  
and show that this structure converges to the strictly commutative 
associative algebra given by the wedge product
on forms. In particular, all of the higher homotopies of the C-infinity
algebra converge to zero.

Only definition ~\ref{defn:cup} and theorem ~\ref{thm:prod} are used in
later sections. 
We begin with the definition of a cochain product on the cochains
of a fixed triangulation $K$.

\begin{defn} \label{defn:cup}
We define $\cup : C^{j}(K) \otimes C^{k}(K) \to C^{j+k}(K)$ by:
\begin{displaymath}
\sigma \cup \tau = R ( W \sigma \wedge W \tau)
\end{displaymath}
\end{defn}

Since $R$ and $W$ are chain maps with respect to $d$ and $\delta$, 
it follows that $\delta$ is a derivation of 
$\cup$, that is, $\delta(\sigma \cup \tau) = \delta \sigma \cup \tau + 
(-1)^{deg (\sigma)}\sigma \cup \delta \tau$. Also, since $\wedge$ is graded 
commutative, $\cup$ is as well: 
$\sigma \cup \tau = (-1)^{deg(\tau)deg(\sigma)} \tau \cup \sigma$. 
It follows from a theorem of Whitney \cite{Wh2} that the product $\cup$
induces the same map on cohomology as the usual (Alexander-Whitney) 
simplicial cochain product. We now give a combinatorial description of $\cup$.

\begin{thm} \label{thm:prod}
Let $\sigma = [p_{\alpha_{0}},p_{\alpha_{1}},\dots,p_{\alpha_{j}}]
\in C^{j}(K)$ and $\tau = [p_{\beta_{0}},p_{\beta_{1}},\dots,p_{\beta_{k}}]
\in C^{k}(K)$. Then $\sigma \cup \tau$ is zero unless $\sigma$ and 
$\tau$ intersect in exactly one vertex and span a $(j+k)$-simplex $\upsilon$,
in which case, for 
$\tau = [p_{\alpha_{j}},p_{\alpha_{j+1}},\dots,p_{\alpha_{j+k}}]$,
we have:
\begin{eqnarray*}
\sigma \cup \tau & = & 
[p_{\alpha_{0}},p_{\alpha_{1}},\dots,p_{\alpha_{j}}] \cup
[p_{\alpha_{j}},p_{\alpha_{j+1}},\dots,p_{\alpha_{j+k}}] \\
& = & s(\sigma, \tau) \frac{j!k!}{(j+k+1)!} [p_{\alpha_{0}},p_{\alpha_{1}},
\dots,p_{\alpha_{j+k}}],
\end{eqnarray*}
where $s(\sigma,\tau)$ is determined by:
\begin{displaymath}
orientation(\sigma) \cdot orientation(\tau) = s(\sigma,\tau) \cdot 
orientation(\upsilon)
\end{displaymath}
\end{thm}

\begin{proof}
Recall that for any simplex $\alpha$,  $W\alpha = 0$ on $M \backslash 
\overline{St(\alpha)}$. So, $\sigma \cup \tau = R(W\sigma \wedge W\tau)$ is
zero if their vertices are disjoint. If $\sigma$ and $\tau$ intersect in more
than one vertex then $W\sigma \wedge W\tau = 0$ since it is a sum of terms 
containing $d \mu_{\alpha_{i}} \wedge d \mu_{\alpha_{i}}$ for some $i$.
Thus, by possibly reordering the vertices of $K$, it suffices to show that 
for $\sigma = [p_{0},p_{1},\dots,p_{j}]$ and  $\tau = [p_{j},p_{j+1},
\dots,p_{j+k}]$, we have that $(\sigma \cup \tau)([p_{0},p_{1},
\dots,p_{j+k}]) = s(\sigma,\tau) \frac{j!k!}{(j+k+1)!}$.
We calculate
{\setlength\arraycolsep{2pt}
\begin{eqnarray*}
\lefteqn{R(W\sigma \wedge W\tau)([p_{0},p_{1},\dots,p_{j+k}])} \\ 
&=&\int_{\upsilon = [p_{0},p_{1},\dots,p_{j+k}]} W([p_{0},p_{1},
\dots,p_{j}]) \wedge W([p_{j},p_{j+1},\dots,p_{j+k}])\\
&=&
j! k! \int_{\upsilon} \
\sum_{i=0}^{j+k} (-1)^{i} \mu_{i} \mu_{j} \ d \mu_{0} \wedge \dots 
\wedge \widehat{d \mu_{i}} \wedge \dots \wedge d \mu_{j+k} \\
\end{eqnarray*}}
Now, $\sum_{i=0}^{j+k} \mu_{i} = 1$, so $d\mu_{0} = - \sum_{i=0}^{j+k} 
d \mu_{i}$, and we have that the last expression
{\setlength\arraycolsep{2pt}
\begin{eqnarray*}
& = & j! k! \int_{\upsilon} \
\sum_{i=0}^{j+k} (-1)^{i} \mu_{i} \mu_{j} \ (- d \mu_{i}) \wedge 
d \mu_{1} \wedge \dots \wedge \widehat{d \mu_{i}} \wedge \dots \wedge 
d \mu_{j+k} \\
& = &  j! k! \int_{\upsilon} \
\mu_{j} \sum_{i=0}^{j+k} \mu_{i} \ d \mu_{1} \wedge \dots \wedge  
d \mu_{j+k} \\
& = &  j! k! \int_{\upsilon}
\mu_{j} \ d \mu_{1} \wedge \dots \wedge  d \mu_{j+k}
\end{eqnarray*}}
Now, $| \int_{\upsilon} d \mu_{1} \wedge \dots \wedge  d \mu_{j+k} |$ 
is the volume of a standard $(j+k)$-simplex, and
thus equals $\frac{1}{(j+k)!}$. From this it is easy to show that 
$\int_{\upsilon} \mu_{j} \ d \mu_{1} \wedge \dots \wedge  d \mu_{j+k} = 
\pm \frac{1}{(j+k+1)!}$, with the appropriate sign
prescribed by the definition of $s(\sigma,\tau)$.
\end{proof}

A special case of this result was derived by Ranicki and Sullivan \cite{RS}
for $K$ a triangulation of a $4k$-manifold and $\sigma$, $\tau$ 
of complimentary dimension. In that paper, they showed
that the pairing given by $\cup$ restricted to simplicies of 
complimentary dimension gives rise to a semi-local combinatorial formula 
for the signature of a $4k$-manifold.

\begin{rmk} \label{rmk:unit}
The constant 0-cochain which evaluates to 1 on all 0-simplicies is the unit 
of the differential graded commutative ring $(C^*, \delta, \cup)$. 
\end{rmk}

We now show that the product $\cup$ converges to $\wedge$, which perhaps
is not surprising, since $\cup$ is induced by the 
Whitney embedding and the wedge product.
Still, the statement may be of computational interest since it shows that
in using cochains to approximate differential forms, the 
product $\cup$ is, in a analytically precise way, an appropriate 
analogue of the wedge product of forms.

\begin{thm} \label{thm:prodconv}
Let $\omega_1,\omega_2 \in \Omega(M)$ and $\sigma$ be an n-simplex of $K$. 
Then there exists a constant $C$ independent of $\omega_1, \omega_2, K$ and
$\sigma$ such that
\[
\lvert W(R\omega_1 \cup R\omega_2)(p) - \omega_1 \wedge \omega_2 (p) \rvert_p
\leq
C \cdot \bigg( c_1 \cdot sup \Big| 
\frac{\partial \omega_2}{\partial x^i} \Big| + c_2 \cdot 
sup \Big| \frac{\partial \omega_1}{\partial x^i} \Big| \bigg) \cdot \eta
\]
for all $p \in \sigma$, where $c_m =  sup \lvert \omega_m \rvert_p$, 
the supremums are taken over all $p \in \sigma$ and $i=1,2, \dots n$ and 
the partial derivatives are taken with respect to a coordinate neighborhood 
containing $\sigma$. 
\end{thm}

\begin{rmk} \label{rmk:prodrmk}
By remark ~\ref{rmk:unit}, theorem ~\ref{thm:prodconv} reduces to
theorem ~\ref{thm:app} in the case $\omega_1$ is the constant function 1. 
\end{rmk}

\begin{proof} Let $\sigma = [p_0 , \dots , p_n]$ be an $n$-simplex 
contained in a coordinate neighborhood with coordinate functions 
$x_1, \dots , x_n$. Let $\mu_i$ denote the $i^{th}$ barycentric coordinate of
$\sigma$.
By the triangle inequality, and a possible reordering of the 
coordinate functions, it suffices to consider the case
{\setlength\arraycolsep{2pt}
\begin{eqnarray*}
\omega_1 &=& f \ d \mu_1 \wedge \dots \wedge d\mu_j \\
\omega_2 &=& g \ d \mu_{\alpha_1} \wedge \dots \wedge d \mu_{\alpha_k} 
\end{eqnarray*}}
We first compute $W(R\omega_1 \cup R\omega_2)$. 
We'll use the notation $[p_s , \dots , p_{s+t}]$ to denote both the simplicial
chain and the simplicial cochain taking the value one on this chain and 
zero elsewhere. Let
{\setlength\arraycolsep{2pt} 
\begin{eqnarray*}
N &=& \{ 0,1,2, \dots ,n \} \\
J &=& \{ 1,2, \dots ,j \} \\
K  &=& \{ \alpha_1, \dots ,\alpha_k \}.
\end{eqnarray*}}
Then
{\setlength\arraycolsep{2pt}
\begin{eqnarray*}
R \omega_1 &=& \sum_{\beta \in N-J}  \Big(
\int_{[p_\beta , p_1, \dots , p_j]} \omega_1 \Big) \ [p_\beta , p_1, \dots , p_j] \\
R \omega_2 &=& \sum_{\gamma \in N-K} \Big( 
\int_{[p_\gamma , p_{\alpha_1}, \dots ,p_{\alpha_k}]} \omega_2 \Big) \
[p_\gamma , p_{\alpha_1}, \dots ,p_{\alpha_k}].
\end{eqnarray*}}
Now, to compute $R \omega_1 \cup R \omega_2$, we use theorem 
~\ref{thm:prod}. If the sets $J$ and $K$ intersect in two or more elements then
$R \omega_1 \cup R \omega_2 = 0$ since, in this case, all products 
of simplicies are zero.

Now suppose that $J$ and $K$ intersect in exactly one element. Without loss
of generality, let us assume $\alpha_1 = 1$. Then the 
product
$$
[p_\beta , p_1, \dots , p_j] \cup 
[p_\gamma , p_{\alpha_1}, \dots ,p_{\alpha_k}]
$$
is non-zero only if $\beta , \gamma$ are distinct elements of the set
$Q =  N - (J \bigcup K)$. Using the abbreviated notation
{\setlength\arraycolsep{2pt}
\begin{eqnarray*}
\ [p_s,p_J,p_K] &=& [p_s , p_1, \dots , p_j, p_{\alpha_1}, \dots ,p_{\alpha_k}]
\\ \int_{[s]} \omega_1 &=& \int_{[p_s , p_1, \dots , p_j]} \omega_1
\\ \int_{[s]} \omega_2 &=& \int_{[p_s , p_{\alpha_1} , \cdots , p_{\alpha_k}]} 
\omega_2 
\end{eqnarray*}}
we compute
$$
R \omega_1 \cup R \omega_2 = \frac{j! k!}{(j+k+1)!}
\sum_{\substack{\beta , \gamma \in Q \\ \beta \neq \gamma}}
\Big(\int_{[\beta]} \omega_1 \Big) \Big( \int_{[\gamma]} \omega_2 \Big) \
[p_\beta, p_\gamma , p_J, p_K].
$$
If all of the coefficients (given by the integrals of $\omega_1$ and $\omega_2$)
were equal, the above expression would vanish, since the terms would cancel
in pairs (by reversing the roles of $\beta$ and $\gamma$). Of course,
this is not the case, but the terms are \emph{almost} equal. We'll
use some estimation techniques developed by Dodziuk and Patodi \cite{DP}.

An essential estimate that we'll need for this case and the next is the
following: there is a constant $c$, independent of $\omega_1, \omega_2, K$ and
$\sigma$, such that for any $p \in \sigma$, and $\beta, \gamma$ as above,
\begin{equation} \label{eq:est1}
\bigg| j!k! \int_{[\beta]} \omega_1 \int_{[\gamma]} \omega_2 - f(p)g(p) \bigg| 
\ \leq \ 
 c \cdot \bigg( c_1 \cdot sup \Big| 
\frac{\partial \omega_2}{\partial x^j} \Big| + c_2 \cdot 
sup \Big| \frac{\partial \omega_1}{\partial x^j} \Big| \bigg) \cdot \eta^{j+k+1}
\end{equation}
where $c_m =  sup \lvert \omega_m \rvert_p$ and the supremums are taken over 
all $p \in \sigma$ and $i=1,2, \dots n$.

To prove this, first note that by the mean value theorem, 
for any points $p,q \in \sigma$,
$\lvert \omega_1(q) - \omega_1(p) \rvert_q \leq c \cdot 
sup |\frac{\partial \omega_1}{\partial x^j} | \cdot \eta$. (Here we're using 
the fact that the Riemannian metric and the flat one induced by pulling back 
along the coordinates $x^i$ are commensurable.) Similarly
for $\omega_2$. Now, fix $p \in \sigma$ and let $dV_\beta$ be the volume 
element on $[p_\beta , p_1, \dots , p_j]$, and $dV_\gamma$ be the volume 
element on $[p_\gamma , p_1, \dots , p_j]$. Then
{\setlength\arraycolsep{2pt}
\begin{eqnarray*} 
\lefteqn{\bigg| j!k! \int_{[\beta]} \omega_1 \int_{[\gamma]} 
\omega_2 - f(p)g(p) \bigg|} \\ 
&=&
\bigg|  j!k! \int_{[\beta]} \omega_1 \int_{[\gamma]} \omega_2 - 
\frac{\int_{[\beta]} f(p) d \mu_1 \wedge \dots \wedge d \mu_j}
{\int_{[\beta]} d \mu_1 \wedge \dots \wedge d \mu_j}
\frac{\int_{[\gamma]} g(p) d \mu_{\alpha_1} \wedge \dots\wedge d \mu_{\alpha_k}}
{\int_{[\gamma]} d \mu_{\alpha_1} \wedge \dots \wedge d \mu_{\alpha_k}} \bigg| \\
&=& j!k! \bigg| \int_{[\beta]} \omega_1 \int_{[\gamma]} \omega_2 -
\int_{[\beta]} \omega_1(p) \int_{[\gamma]} \omega_2(p) \bigg| \\
&\leq& j!k! \bigg| \int_{[\beta]} \omega_1 \bigg|
 \bigg| \int_{[\gamma]} \omega_2 - \int_{[\gamma]} \omega_2(p) \bigg|
+ \bigg| \int_{[\gamma]} \omega_2 \bigg|
 \bigg|\int_{[\beta]} \omega_1 - \int_{[\beta]} \omega_1(p) \bigg| \\
&\leq& j!k! \ c_1 \cdot \eta^j \int_{[\gamma]} \lvert \omega_2 - \omega_2(p) 
\rvert_q \ dV_{\gamma} + 
c_2 \cdot \eta^k \int_{[\beta]} \lvert \omega_1 
- \omega_1(p) \rvert_q \ dV_{\beta}\\
&\leq& c \cdot \bigg( c_1 \cdot sup \Big| 
\frac{\partial \omega_2}{\partial x^i} \Big| +  c_2 
\cdot sup \Big| \frac{\partial \omega_1}{\partial x^i} \Big| \bigg) 
\cdot \eta^{j+k+1}.
\end{eqnarray*}}
This implies, by the triangle inequality, for any $\beta,\gamma$
\begin{equation}  \label{eq:coefest}
\bigg| \int_{[\beta]} \omega_1 \int_{[\gamma]} \omega_2 - 
 \int_{[\gamma]} \omega_1 \int_{[\beta]} \omega_2
\bigg|
\ \leq \ 
 c \cdot \bigg( c_1 \cdot sup \Big| 
\frac{\partial \omega_2}{\partial x^j} \Big| + c_2 \cdot 
sup \Big| \frac{\partial \omega_1}{\partial x^j} \Big| \bigg) \cdot \eta^{j+k+1}
\end{equation}

Now that we have estimated the coefficients of $W( R \omega_1 \cup R\omega_2)$,
this case is completed by estimating the the product of the $d \mu_i$'s that
appear in $W( R \omega_1 \cup \omega_2)$. As shown in \cite{DP},
$$
| d \mu_i |_p \leq \frac{\lambda}{r ( p_i , | \sigma_i| )}, 
$$
where $\sigma_i = [p_0, \cdots, p_{j-1} , p_{j+1}, \cdots , p_N]$ is the face 
opposite of $p$, and $r$ is the Riemannian geodesic distance. So, by lemma
\ref{lemma:techass}
$$
| d \mu_i |_p \leq \lambda' \cdot \eta^{-1}
$$
for some constant $\lambda'$, and therefore 
\begin{equation} \label{eq:dmuest}
| d \mu_{i_1} \wedge \cdots \wedge d \mu_{i_{j+k}} |_p \leq
| d \mu_{i_1} |_p \dots | d \mu_{i_{j+k}} |_p \leq \lambda \cdot \eta^{-(j+k)}.
\end{equation}
By combining \eqref{eq:coefest} and \eqref{eq:dmuest}, we finally have, 
for the case that $J$ and $K$ intersect in exactly one element, 
{\setlength\arraycolsep{2pt}
\begin{eqnarray*}
| W( R \omega_1 \cup R \omega_2)(p) - \omega_1 \wedge \omega_2 (p) |_p &=&
 | W( R \omega_1 \cup R \omega_2)(p)|_p \\ &\leq& 
C \cdot \bigg( c_1 \cdot sup \Big| 
\frac{\partial \omega_2}{\partial x^i} \Big| + c_2 \cdot 
sup \Big| \frac{\partial \omega_1}{\partial x^i} \Big| \bigg) \cdot \eta
\end{eqnarray*}}

We now consider the case that $J$ and $K$ are disjoint. 
We first note that for any $\tau \in Q = N - (J \cup K)$, 
there are exactly $j+k+1$ products
$$
[p_\beta, p_1, \dots , p_j] \cup [ p_\gamma, p_{\alpha_1}, \dots ,p_{\alpha_k}]
$$
which equal a nonzero multiple of $[p_\tau , p_J , p_K] = 
[p_\tau, p_1, \dots , p_j, p_{\alpha_1}, \dots ,p_{\alpha_k}]$. 
These are given by the three mutually exclusive cases:
{\setlength\arraycolsep{2pt}
\begin{eqnarray*}
\beta &=& \tau, \  \gamma \in J \\
\gamma &=& \tau, \  \beta \in K \\
\beta &=& \gamma = \tau
\end{eqnarray*}}
Using the same notation as the previous case, we compute 
{\setlength\arraycolsep{2pt}
\begin{eqnarray} \label{eq:thingy}
R \omega_1 \cup R \omega_2 &=& \frac{j! k!}{(j+k+1)!} \bigg( \
\sum_{\|0\|} \Big(\int_{[\beta]} \omega_1 \Big) 
\Big( \int_{[\gamma]} \omega_2 \Big) \ [p_0, p_J, p_K] \\
&&
{}+\sum_{\tau \in Q-\{0\}} \sum_{\|\tau\|} \Big(\int_{[\beta]} \omega_1 \Big) 
\Big( \int_{[\gamma]} \omega_2 \Big) \ [p_\tau, p_J,p_K] \bigg) \nonumber
\end{eqnarray}}
where the sums labeled $\displaystyle \sum_{\|s\|}$ are over all 
$\beta, \gamma$ such that
$$
[p_\beta, p_1, \dots , p_j] \cup [ p_\gamma, p_{\alpha_1}, \dots ,p_{\alpha_k}]
= \frac{j! k!}{(j+k+1)!} [p_s, p_J , p_K].
$$ 
From lemma \ref{lemma:Wcomp}, which follows the proof of this 
theorem, 
\begin{eqnarray*}
W([p_0, p_J, p_K]) = (j+k)! \ d\mu_J \wedge d\mu_K - \sum_{r \in Q -\{0\}} 
W([p_\tau, p_J , p_K])
\end{eqnarray*}
So, 
{\setlength\arraycolsep{2pt}
\begin{eqnarray*}
\lefteqn{| W(R \omega_1 \cup R \omega_2) - \omega_1 \wedge \omega_2 |_p} \\
&\leq& 
\frac{j! k!}{(j+k+1)} \ \bigg|
\sum_{\|0\|} \Big( \int_{[\beta]} \omega_1 \Big) \Big( \int_{[\gamma]} \omega_2 
\Big) \ d\mu_J \wedge d\mu_K - \omega_1 \wedge \omega_2 \bigg|_p \\
&& + \frac{j!k!}{(j+k+1)!} \
\bigg|\sum_{\tau \in Q-\{0\}} \bigg( \sum_{\|\tau\|} 
\Big(\int_{[\beta]} \omega_1 \Big) 
\Big( \int_{[\gamma]} \omega_2 \Big) \\ && \hspace{130pt}
-\sum_{\|0\|} \Big(\int_{[\beta]} \omega_1 \Big) 
\Big( \int_{[\gamma]} \omega_2 \Big) \bigg) \ W([p_\tau,p_J, p_K]) \bigg|_p
\end{eqnarray*}}
By our estimates in \eqref{eq:coefest} and \eqref{eq:dmuest}, 
the latter term is bounded appropriately.
As for the first term, recall that the sum $\displaystyle \sum_{\|0\|} $ 
consists of j+k+1 terms. We use \eqref{eq:coefest} again to bound
\begin{equation} \label{eq:est2}
\bigg| \sum_{\|0\|} \Big( \int_{[\beta]} \omega_1 \Big) \Big( \int_{[\gamma]} 
\omega_2 \Big) 
- (j+k+1)\Big( \int_{[0]} \omega_1 \Big) \Big( \int_{[0]} \omega_2 
\Big) \bigg|
\end{equation}
and using \eqref{eq:est1}, for fixed $p \in \sigma$ we have a bound on
\begin{equation} \label{eq:est3}
\bigg| \Big( \int_{[0]} \omega_1 \Big) \Big( \int_{[0]} \omega_2 
\Big) - f(p)g(p) \bigg|.
\end{equation}
Finally, using the triangle inequality and 
combining \eqref{eq:est2} and \eqref{eq:est3} with \eqref{eq:dmuest}
we can conclude
{\setlength\arraycolsep{2pt}
\begin{eqnarray*}
&& \hspace{-100pt} \frac{j! k!}{(j+k+1)} \ \bigg|
\sum_{\|0\|} \Big( \int_{[\beta]} \omega_1 \Big) \Big( \int_{[\gamma]} \omega_2 
\Big) \ d\mu_J \wedge d\mu_K (p)- \omega_1 \wedge \omega_2 (p) \bigg|_p
\\ &\leq&
C \cdot \bigg( c_1 \cdot sup \Big| 
\frac{\partial \omega_2}{\partial x^j} \Big| + c_2 \cdot 
sup \Big| \frac{\partial \omega_1}{\partial x^j} \Big| \bigg) \cdot \eta
\end{eqnarray*}}
\end{proof}

\begin{lemma} \label{lemma:Wcomp} Let $\sigma = [p_0 , p_1 , \dots ,p_n]$,
$N = \{1,2,\dots,n\}$ and $I = \{ i_1 , \dots i_m\} \subset N$. Then
$$
W([p_0, p_{i_1}, \dots , p_{i_m}]) =  m! \ d \mu_{i_1} \wedge \cdots 
\wedge d\mu_{i_m} - \sum_{r \in N - I} W([p_r ,  p_{i_1}, \dots , p_{i_m}])
$$
\end{lemma}

\begin{proof} The proof is a computation. We let
{\setlength\arraycolsep{2pt}
\begin{eqnarray*}
d\mu_I &=& d\mu_{i_1} \wedge \cdots \wedge d\mu_{i_m} \\
d\mu_I^s &=& d\mu_{i_1} \wedge \cdots \wedge \widehat{d\mu_{i_s}} \wedge \cdots 
\wedge d\mu_{i_m}
\end{eqnarray*}}
and compute
{\setlength\arraycolsep{2pt}
\begin{eqnarray*}
\frac{1}{m!} W([p_0, p_{i_1}, \dots , p_{i_m}]) &=&
\mu_0 \ d\mu_I + \sum_{s=1}^m (-1)^s \mu_{i_s} \ d\mu_0 \wedge d\mu^s_I \\
&=& \Big(1-\sum_{r=1}^n \mu_r \Big) \ d\mu_I + 
\sum_{s=1}^m (-1)^s \mu_{i_s} 
\ \Big( - \sum_{r=1}^n d \mu_r  \Big) \wedge d\mu^s_I\\
&=& d\mu_I - \sum_{r=1}^n \mu_r \ d\mu_I - \sum_{s=1}^m (-1)^s \mu_{i_s} \ 
\Big(d\mu_{i_s} + \sum_{r \in N - I} d\mu_r \Big) \wedge d\mu_I^s \\
&=&
d\mu_I - \sum_{r \in N - I} \mu_r \ d\mu_I - \sum_{s=1}^m (-1)^s \mu_{i_s} 
\Big( \sum_{r \in N - I} d\mu_r \Big) \wedge d\mu_I^s \\
&=& d\mu_i - \sum_{r \in N - I} \Big(\mu_r \ d\mu_I + \sum_{s=1}^m (-1)^s 
\mu_{i_s} \ d\mu_r \wedge d\mu_I^s \Big) \\
&=& d\mu_I - \frac{1}{m!} \sum_{r \in N - I} 
W([p_r ,  p_{i_1}, \dots , p_{i_m}])
\end{eqnarray*}}
\end{proof}

\begin{cor} \label{cor:prodconv2}
There exists a constant $C$ and positive integer $m$, independent of 
$K$ such that
\[
\| W(R\omega_1 \cup R\omega_2) - \omega_1 \wedge \omega_2 \|
\leq
C \cdot \lambda ( \omega_1,\omega_2 ) \cdot \eta
\]
where
\[
\lambda (\omega_1 , \omega_2 ) = \|\omega_1\|_{\infty} \cdot 
\|(Id + \Delta)^m \omega_2 \| +  \|\omega_2\|_{\infty} \cdot
\|(Id + \Delta)^m \omega_1 \| 
\]
for all smooth forms $\omega_1,\omega_2 \in \Omega(M)$, where 
$\| \hspace{1em} \|$ is the $\mathcal{L}_2$-norm on $M$.
\end{cor}

\begin{proof} We integrate the point-wise estimate from theorem 
\ref{thm:prodconv}, using the facts that $M$ is compact, $sup \lvert \omega_k
\rvert = \| \omega_k \|_{\infty}$, and the Sobolev-Inequality
\[
sup \Big| \frac{\partial \omega_k}{\partial x^i} \Big|
\leq C \cdot \| \omega_k \|_{2m} = C \cdot \| (Id + \Delta)^m \omega_k \|
\]
for sufficiently large $m$, where $\| \hspace{1em} \|_{2m}$ is the Sobolev
$2m$-norm.
\end{proof}

The convergence of $\cup$ to the associative product $\wedge$ is, a priori, 
a bit mysterious due to the following:

\begin{ex} The product $\cup$ is not associative. For example, in 
the figure below, $(a \cup b) \cup e = 0$, since $a$ and $b$ 
do not span a $0$-simplex, but $a \cup (b \cup e) = -\frac{1}{4} e$.
\end{ex}

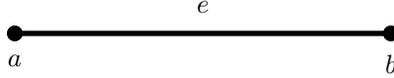
\begin{center}
\begin{pspicture}(0,0)(5,2)
\psline[linewidth=2pt]{*-*}(0,1)(5,1)
\psset{labelsep=8pt}
\uput[d](0,1){$a$}
\uput[u](2.5,1){$e$}
\uput[d](5,1){$b$}
\end{pspicture}
\end{center}

In the above example, the cochains $a,b$ and $e$ may be thought of as
delta functions, in the sense that they evaluate to one on a single simplex
and zero elsewhere. If we work with cochains which are ``smoother'', i.e.
represented by the integral of a smooth differential form,
associativity is \emph{almost} obtained. In fact, the next theorem shows that 
for such cochains, the deviation from being associative is
bounded by a constant times the mesh of the triangulation. Hence, associativity
is recovered in the mesh goes to zero limit. 

\begin{thm} \label{thm:nearassoc}
There exists a constant $C$ and positive integer $m$, independent of 
$K$ such that
\[
\| (R\omega_1 \cup R\omega_2) \cup R\omega_3 - R\omega_1 (R\omega_2 \cup 
R\omega_3) \| \leq
C \cdot \lambda ( \omega_1,\omega_2,\omega_3 ) \cdot \eta
\]
for all $\omega_1 , \omega_2,\omega_3 \in \Omega(M)$, 
where $\| \hspace{1em} \|$ is the Whitney norm and
\[
\lambda (\omega_1 , \omega_2,\omega_3 ) = 
\sum \|\omega_r\|_{\infty} \cdot \|\omega_s\|_{\infty} \cdot
\|(Id + \Delta)^m \omega_t \|
\]
where the sum is over all cyclic permutations $\{r,s,t\}$ of $\{1,2,3\}$.
\end{thm}

\begin{proof} We can prove this by first showing each of 
$(R\omega_1 \cup R\omega_2) \cup R\omega_3$ and $R\omega_1 \cup (R\omega_2 \cup 
R\omega_3)$ are close to $\omega_1 \wedge \omega_2 \wedge \omega_3$ 
in the point-wise norm $\lvert \hspace{1em} \rvert_p$. The final result
is then obtained by integrating and applying the Sobolev inequality to each
point-wise error, then applying the triangle inequality.

Let $A \approx B$ mean 
$$ \lvert A - B \rvert_p \leq c \cdot \sum \| \omega_r \|_{\infty} \cdot 
\| \omega_s \|_{\infty} \cdot sup \Big| 
\frac{\partial \omega_t}{\partial x^i} \Big| \cdot \eta
$$
We'll consider the first case, 
{\setlength\arraycolsep{2pt}
\begin{eqnarray} \label{eq:showthis}
W ((R\omega_1 \cup R\omega_2) \cup R\omega_3) \approx \omega_1 \wedge \omega_2, 
\end{eqnarray}}
only; the second case is similar.

It suffices to consider the case
{\setlength\arraycolsep{2pt}
\begin{eqnarray*}
\omega_1 &=& f \ d \mu_1 \wedge \dots \wedge d\mu_j \\
\omega_2 &=& g \ d \mu_{\alpha_1} \wedge \dots \wedge d \mu_{\alpha_k} \\
\omega_3 &=& h \ d \mu_{\beta_1} \wedge \dots \wedge d \mu_{\beta_l}.
\end{eqnarray*}}
The proof is analogous to that of theorem \ref{thm:prodconv},
the only differences are that the combinatorics of two cochain products
is slightly more complicated, and the estimates now involve coefficients
which are triple products of integrals over simplicies. Let
{\setlength\arraycolsep{2pt}
\begin{eqnarray*}
N &=& \{1, \dots ,n\} \\
J &=& \{1, \dots ,j\} \\
K &=& \{ \alpha_1, \dots ,\alpha_k \} \\
L &=& \{ \beta_1, \dots ,\beta_l \} \\
Q &=& N - (J \cup K \cup L)
\end{eqnarray*}}
Let us assume $J \cap K \cap L = \emptyset$; the other cases are similar.
Define $A \sim B$ by 
$$
|A - B| 
\leq C \cdot \|\omega_r\|_{\infty} \cdot \|\omega_s\|_{\infty} 
\cdot sup \Big| 
\frac{\partial \omega_t}{\partial x^i} \Big| \cdot \eta^{j+k+l+1}
$$
Using similar techniques as in the proof of theorem ~\ref{thm:prodconv},
for all $a \in N - J$, $b \in N - K$, $c \in N - L$
{\setlength\arraycolsep{2pt}
\begin{eqnarray} \label{eq:tripcoef}
\Big(\int_{[a]} \omega_1 \Big) \Big( \int_{[b]} \omega_2 \Big) 
\Big(\int_{[c]} \omega_3 \Big) &\sim&
\Big(\int_{[0]} \omega_1 \Big) \Big( \int_{[0]} \omega_2 \Big) 
\Big(\int_{[0]} \omega_3 \Big) \\ \nonumber
j!k!l! \Big(\int_{[0]} \omega_1 \Big) \Big( \int_{[0]} \omega_2 \Big) 
\Big( \int_{[0]} \omega_3 \Big) &\sim& f(p)g(p)h(p)
\end{eqnarray}}
For any $\tau \in Q$, there are exactly
$$
(j+k+1)(j+k+1) + (j+k+1)l = (j+k+1)(j+k+l+1)
$$
products
$$
[p_a, p_1, \dots , p_j] \cup [ p_b, p_{\alpha_1}, \dots ,p_{\alpha_k}]
\cup [ p_c, p_{\beta_1}, \dots ,p_{\beta_l}]
$$
that equal a non-zero multiple of $[ p_\tau, p_J, p_K, p_L]$.
Then
$$
\frac{j! k! (j+k!) l!}{(j+k+1)!(j+k+l+1)!} (j+k+1)(j+k+l+1) 
= \frac{j! k! l!}{(j+k+l)!}
$$
so that, by applying lemma \ref{lemma:Wcomp},
and equations \eqref{eq:tripcoef} and \eqref{eq:dmuest},
{\setlength\arraycolsep{2pt}
\begin{eqnarray*}
\lefteqn{W((R \omega_1 \cup R \omega_2) \cup R\omega_3)} \\
&\approx&
\frac{j! k! l!}{(j+k+l)!} 
\bigg( \Big(\int_{[0]} \omega_1 \Big) \Big( \int_{[0]} \omega_2 \Big) 
\Big( \int_{[0]} \omega_3 \Big) \ W([p_0, p_J, p_K,p_L]) \\
&&
{}+ \sum_{\tau \in Q-\{0\}}
\Big(\int_{[0]} \omega_1 \Big) \Big( \int_{[0]} \omega_2 \Big) 
\Big( \int_{[0]} \omega_3 \Big) \ 
W([p_\tau, p_J,p_K,p_L]) \bigg) \nonumber
\end{eqnarray*}}

\end{proof}

In the previous theorem, we dealt with the non-associativity of $\cup$ 
analytically. There is also an algebraic way to deal with this, via
an algebraic generalization of commutative, associative algebras, called
$\mathcal{C}_\infty$-algebras. First we'll give an abstract definition,
and then unravel what it means.

\begin{defn} Let $C$ be a graded vector space, and $C[-1]$ denote the
graded vector space $C$ with grading shifted down by one. Let 
$\mathcal{L}(C) = \bigoplus_i \mathcal{L}^i(C)$ be the 
free Lie co-algebra on C. 
A $\mathcal{C}_\infty$-algebra on $C$ is a degree 1 co-derivation 
$D:\mathcal{L}(C[-1]) \to \mathcal{L}(C[-1])$ such that $D^2 = 0$.
\end{defn}

A co-derivation on a free Lie co-algebra is uniquely determined by
a collection of maps from $\mathcal{L}^i(C)$ to $C$ 
for each $i \geq 1$.
If we let $m_i$ denote the restriction of $D$ to $\mathcal{L}^i(C)$, then
the equation $D^2$ = 0 is equivalent to a collection of equations:
{\setlength\arraycolsep{2pt}
\begin{eqnarray*}
m_1^2 &=& 0 \\
m_1 \circ m_2 &=& m_2 \circ m_1 \\
m_2 \circ m_2 - m_2 \circ m_2 &=& m_1\circ m_3 + m_3 \circ m_1 \\ 
&\vdots&
\end{eqnarray*}}
We can regard $m_1$ as a differential and $m_2$ a commutative 
multiplication on $C$. The second equation states that $m_1$ is a
derivation of $m_2$. The third equation states that $m_2$ is 
associative up to the (co)-chain homotopy $m_3$. Note that, due to the 
shift of grading, $m_j$ has degree $2-j$.

The following theorem is due to Sullivan \cite{DS}. See also \cite{TT}
for use of similar techniques.

\begin{thm} \label{thm:sully}
Let $(C,\delta)$ be the simplicial cochains of a triangulated space and
$\cup$ be any local commutative (possibly non-associative) 
cochain multiplication
on $C$ such that $\delta$ is a derivation of $\cup$. Then there is a 
canonical local inductive construction which extends $(C,\delta,\cup$)
to a $\mathcal{C}_\infty$-algebra.
\end{thm}

In the theorem, local means that the product of a $j$-simplex and a $k$-simplex
is zero unless they span a $j+k$-simplex, in which case it is a multiple
of this simplex. By theorem ~\ref{thm:prod}, the commutative
product $\cup$ defined at the beginning of this section satisfies this
and the other conditions of theorem ~\ref{thm:sully}. 

The next theorem shows that the $\mathcal{C}_\infty$-algebra on $C$ converges 
to the strict commutative and associative algebra given by 
the wedge product on forms in a sense analogous to the convergence statements 
we've made previously. In particular, all higher homotopies converge to zero
as the mesh tends to zero.

\begin{thm} \label{thm:hmtps2zero}
Let $C$ be the simplicial cochains of a triangulation $K$ of $M$, 
with mesh $0 \leq \eta \leq 1$. Let $m_1 = \delta , m_2 = \cup ,
m_3 , \dots$ be the extension of $C,\delta ,\cup$  to a $\mathcal{C}_\infty$-
algebra as in theorem ~\ref{thm:sully}.
Then there exists a constant $\lambda$ independent of $K$ such that,
for all $j \geq 3$,
$$
\| W( m_j ( R\omega_1 , \dots , R\omega_j)) \| \leq
\lambda \cdot \prod_{i=1}^j \|\omega_i\|_{\infty} \cdot \eta
$$
for all $\omega_1, \dots ,\omega_k \in \Omega(M)$.
\end{thm}

\begin{proof} 
Suppose degree $\omega_1, \dots ,\omega_j$ are of degree
$\alpha_1 , \dots , \alpha_j$, respectively. Let $\alpha = \sum \alpha_i$.
We need two facts. First, for any $\alpha_i$-simplex $\tau$ of $K$,
\begin{equation} \label{eq:first}
| R \omega_i (\tau) | \leq c \cdot \|\omega_i\|_{\infty} \cdot \eta^{\alpha_i}
\end{equation}
Secondly, if $p$ is a point in an $n$-simplex $\sigma$, and the $r$-simplicies
which are faces of $\sigma$ are $\sigma_{r}^1 , \dots , \sigma_{r}^m$
then, by equation \eqref{eq:dmuest}, 
\begin{equation} \label{eq:two}
\bigg| W \bigg( \sum_{i=1}^m \sigma_{r}^i \bigg) \bigg|_p 
\leq c' \cdot \eta^{-r}.
\end{equation}
Now, since $m_j$ has degree $2-j$, $m_j( R\omega_1 , \dots , R\omega_j)$
is a linear combination of $(\alpha + 2 - j)$-simplicies. 
Combining this with \eqref{eq:first} and \eqref{eq:two}, we have for all 
$p\in M$ and some $\lambda \geq 0$
{\setlength\arraycolsep{2pt}
\begin{eqnarray*}
| W( m_j ( R\omega_1 , \dots , R\omega_j)) |_p &\leq& 
\lambda \cdot \prod_{i=1}^j  \|\omega_i\|_{\infty} \cdot \eta^\alpha \cdot 
\eta^{-(\alpha + 2 - j)} \\
&\leq& \lambda \cdot \prod_{i=1}^j  \|\omega_i\|_{\infty} \cdot \eta
\end{eqnarray*}}
The result is obtained by integrating over $M$.
\end{proof}

\section{Combinatorial Star Operator} \label{sec:star}

In this section we define the combinatorial star operator $\bigstar$ and prove
that it provides a good approximation to the smooth Hodge-star $\star$.
We also examine the relations which are expected to hold by analogy with the
smooth setting. We find that some hold precisely,
while others may only recovered in the mesh goes to zero limit.

\begin{defn} \label{defn:bigstar}
Let $K$ be a triangulation of a closed orientable manifold $M$,
with simplicial cochains $C = \bigoplus_j C^j$. 
Let $\langle , \rangle$ be a non-degenerate positive definite 
inner product on $C$ such that
$C^i \perp C^j$ for $i \neq j$. 
For $\sigma \in C^{j}$ we define $\bigstar \sigma \in C^{n-j}$ by:
\begin{displaymath}
\langle \bigstar \sigma, \tau \rangle = (\sigma \cup \tau)[M]
\end{displaymath}
where $[M]$ denotes the fundamental class of $M$.
\end{defn}

We emphasize that, as motivated by definition~\ref{defn:star}, the
essential ingredients of a star operator are Poincar\'e Duality and a
non-degenerate inner product. We can regard the inner product as giving
some geometric structure to the space. In particular it gives lengths of,
and angles between, edges. As in the smooth setting, the star operator
depends on the choice of inner product (or Riemannian metric).
See section~\ref{sec:ip} for the definition of a particularly nice class 
of inner products that we call \emph{geometric inner products}.

Here are some elementary properties of $\bigstar$. 

\begin{thm} \label{thm:starprop} The following hold:
\begin{enumerate}
\item $\bigstar \delta = (-1)^{j+1} \delta^{\ast} \bigstar$, \ i.e. $\bigstar$
is a chain map.
\item  For $\sigma \in C^j$ and $\tau \in C^{n-j}$,
$\langle \bigstar \sigma , \tau \rangle 
= (-1)^{j(n - j)} \langle \sigma, \bigstar \tau \rangle$, i.e. 
$\bigstar$ is (graded) skew-adjoint.
\item $\bigstar$ induces isomorphisms $\mathcal{H}C^{j}(K) \to 
\mathcal{H}C^{n-j}(K)$ on harmonic cochains.
\end{enumerate}
\end{thm}

\begin{proof} The first two proofs are computational:
\begin{enumerate}
\item For $\sigma, \tau \in C$, we have:
{\setlength\arraycolsep{2pt}
\begin{eqnarray*}
\langle \bigstar \delta \sigma, \tau \rangle 
& = & (\delta \sigma \cup \tau)[M] \\
& = & (-1)^{j+1}(\sigma \cup \delta \tau)[M] \\ 
& = & (-1)^{j+1} \langle \bigstar \sigma, \delta \tau \rangle \\
& = & \langle (-1)^{j+1} \delta^{\ast} \bigstar \sigma, \tau \rangle
\end{eqnarray*}}
where we have used that fact that $d$ is a derivation of $\cup$ and 
$M$ is closed.
\item We compute:
{\setlength\arraycolsep{2pt}
\begin{eqnarray*}
\langle \bigstar \sigma, \tau \rangle 
& = & (\sigma \cup \tau)[M] \\
& = & (-1)^{j(n-j)}(\tau \cup \sigma)[M] \\ 
& = & (-1)^{j(n-j)} \langle \bigstar \tau, \sigma \rangle \\
& = & (-1)^{j(n-j)} \langle \sigma, \bigstar \tau \rangle
\end{eqnarray*}}
\item Via the Hodge-decomposition of cochains, $\mathcal{H}C^{j}(K)$ may
be identified with the cohomology $\mathcal{H}^{j}(K)$. Here $\bigstar$
is the composition of two isomorphisms, Poincare Duality (since $M$ is a 
manifold) and the inverse of the non-degenerate metric.
\end{enumerate}
\end{proof}

We remark here that $\bigstar$ is in general not invertible, since the 
cochain product does not necessarily give rise to a non-degenerate pairing 
(on the cochain level!). This implies that $\bigstar$ is not an orthogonal
map, and $\bigstar^2 \neq \pm Id$. 

For the remainder of this section, we'll fix the inner product on
cochains to be the Whitney inner product, so that $\bigstar$ is the 
star operator induced by the Whitney inner product.
This will
be essential in showing that $\bigstar$ converges to the smooth Hodge star 
$\star$, which is defined using the Riemannian metric. First, a useful lemma. 
Let $\perp$ denote the orthogonal projection of $\Omega^{j}(M)$ onto 
the image of $C^{j}(K)$ under the Whitney embedding $W$.

\begin{lemma} $W \bigstar = \perp \star W$
\end{lemma}

\begin{proof} 
Let $a \in C^{j}(K)$ and $b \in C^{n-j}(K)$. Note that $\star W a$ is an
$\mathcal{L}_{2}$-form but in general is not a Whitney form. We compute:
\begin{displaymath}
\langle W \bigstar a , W b \rangle = \langle \bigstar a, b \rangle
= \int_{M} Wa \wedge Wb = \langle \star W a , W b \rangle, 
\end{displaymath}
Thus, $W \bigstar a$ and $\star W a$ have the same inner product with all 
forms in the image of $W$, so $W \bigstar = \perp \star W$.
\end{proof}

Now for our convergence theorem of $\bigstar$:

\begin{thm} \label{thm:starapprox} Let $M$ be a Riemannian manifold with
triangulation $K$ of mesh $\eta$.
There exists a positive constant $C$ and a positive
integer $m$, independent of $K$, such that
$$
\lVert \star \omega - W \bigstar R \omega \rVert \leq C \cdot
\| (Id + \Delta)^m \omega \| \cdot \eta  
$$
for all $C^{\infty}$ differential forms $\omega$ on $M$.
\end{thm}

\begin{proof}
We compute and use Theorem ~\ref{thm:app}
\begin{eqnarray*}
\lVert \star \omega - W \bigstar R \omega \rVert & = & 
\lVert \star \omega - \perp \star W R \omega \rVert \\ 
& \le & \lVert \star \omega - \star W R \omega \rVert
+ \lVert \star W R \omega - \perp \star W R \omega \rVert \\
& \le & \lVert \star \rVert \lVert \omega - W R \omega \rVert
+ \lVert \star W R \omega - W R \star \omega \rVert \\
& \le & \lVert \omega - W R \omega \rVert +
\lVert \star W R \omega - \star \omega \rVert +
\lVert \star \omega - W R \star \omega \rVert \\
& \le & 2 \lVert \omega - W R \omega \rVert
+ \lVert \star \omega - W R \star \omega \rVert \\
& \le & 3 C \cdot \| (Id + \Delta)^m \omega \| \cdot \eta
\end{eqnarray*}
\end{proof}

The operator $\bigstar$ also respects the Hodge decompositions of $C(K)$ and 
$\Omega(M)$ in the following sense:

\begin{thm} \label{thm:happstar}
Let $M$ be a Riemannian manifold with
triangulation $K$ of mesh $\eta$. 
Let $\omega \in \Omega^{j}(M)$, $R\omega \in C^{j}(K)$ have Hodge 
decompositions
\begin{eqnarray*}
\omega & = & d \omega_{1} + \omega_{2} + d^{\ast} \omega_{3}  \\ 
R\omega & = & \delta a_{1} + a_2 + \delta^{\ast} a_3
\end{eqnarray*}
There exists a positive constant $C$ and a positive
integer $m$, independent of $K$, such that
\begin{eqnarray*}
\| \star \omega_{1} - W \bigstar a_{1} \| & \leq & C \cdot 
\| (Id + \Delta)^{m} \omega\| \cdot \eta \\
\| \star d\omega_{2} - W \bigstar \delta a_{2} \| 
& \leq & C \cdot \| (Id + \Delta)^{m} \omega\| \cdot \eta \\
\| \star d^{\ast} \omega_{3} - W \bigstar \delta^{\ast} a_{3} \| 
& \leq & C \cdot \| (Id + \Delta)^{m} \omega \| \cdot \eta 
\end{eqnarray*}
\end{thm}

\begin{proof} The proof is analogous to the proof of 
theorem~\ref{thm:starapprox}.
\end{proof}

One might now ask further questions about convergence, say for 
compositions of the operators $\delta$, $\delta^{\ast}$, and $\bigstar$. 
We now discuss some of these questions.

We first note $\delta$ provides a good approximation of $d$ 
in the sense that $\|d \omega - W \delta R \omega \|$
is bounded by a constant times the mesh. This follows immediately from 
theorem ~\ref{thm:app}, using the fact that $\delta R = R d$. 
In the same way, using theorem ~\ref{thm:starapprox}, $\bigstar \delta$ 
provides a good approximation to $\star d$. In summary, we have:
$$
\pm \delta^{\ast} \bigstar = \bigstar \delta \to 
\star d = \pm d^{\ast} \star
$$

One would also like to know if either of $\delta \bigstar$
or $\bigstar \delta^{\ast}$ provide a good approximation of 
$\pm d^{\ast} \star = d \star$. Answers to these questions are seemingly 
harder to come by. 

As a precursor, we point out that there is not a complete answer 
as to whether or not 
$\delta^{\ast}$ converges to $d^{\ast}$. In \cite{SL}, Smits does prove 
convergence for the case of $1$-cochains on a surface. To the 
author's mind, and as can be seen in the work of \cite{SL},
one difficulty (with the general case) is that the operator $\delta^{\ast}$
is not local, since it involves the inverse of the cochain inner 
product.\footnote{If the cochain inner product is written as a matrix M with
respect to the basis given by the simplicies, then 
$\delta^{\ast} = M^{-1} \partial M$ where $\partial$ is the usual
boundary operator on chains.} A first attempt to understand this
inverse is described in section~\ref{sec:ip}.

The issue becomes further complicated when considering the operator 
$\bigstar \delta^{\ast}$. We have no convergence statements
about this operator. On the other hand, 
the operator $\delta \bigstar$, which incidentally does
not equal $\pm \bigstar \delta^{\ast} $, is a bit less mysterious,
and we have weak convergence in the sense that
$$
\langle W \delta \bigstar R \omega_1 - d \star \omega_1 , \omega_2 \rangle
$$
is bounded by a constant $\lambda$ (depending on $\omega_1$ and 
$\omega_2$) times the mesh.

Finally, one might ask if $\bigstar^2$ approaches $\pm Id$ for a fine
triangulation. While we have no analytic result to state, our 
calculations for the circle in section~\ref{sec:circle} 
suggest this is the case. 
One can show that a graded symmetric operator squares to $\pm Id$ if and only
if it is orthogonal. Hence one might view $\bigstar^2 \neq Id$ as the 
failure of orthogonality, which at least for applications  
to surfaces in section~\ref{sec:surfaces}, presents no difficulty.

\section{Applications to Surfaces} \label{sec:surfaces}

In this section we study applications of the combinatorial star operator 
on a triangulated closed surface. As motivation, let us first recall some 
facts from the analytic setting.

Let $M$ be a Riemann surface. There is a Hodge-star operator in the complex 
valued $1$-forms of $M$, defined in local coordinates by
$\star dx = dy$ and $\star dy = - dx$ and extended linearly over $\mathbb{C}$.
One can check that this is well defined using the Cauchy-Riemann equations for 
the coordinate interchanges. The Hodge-star operator restricts to an 
orthogonal automorphism of complex valued $1$-forms that squares to $-Id$. 
Furthermore, the harmonic $1$-forms split into an orthogonal sum of 
holomorphic and anti-holomorphic $1$-forms corresponding to the $-i$ and 
$+i$ eigenspaces of the Hodge-star operator. 

Riemann studied how the integrals of holomorphic and anti-holomorphic
$1$-forms, called periods, are related to the underlying complex structure.
He showed that for any fixed homology basis these periods satisfy the 
so-called \emph{bi-linear relations}. Furthermore, choosing
a particular basis for the holomorphic $1$-forms gives rise to a 
\emph{period matrix}, which, by Torelli's theorem, 
determines the conformal structure of the Riemann surface. These period
matrices lie in what is called the Siegel upper half space.
(Two references for this material are \cite{GS} and \cite{FK}.)
An unsolved problem, called the Schottky problem, 
is to determine which points in the Siegel upper half space represent 
the period matrix of a Riemann surface.

In this section we'll show that the combinatorial Hodge-star operator
on a triangulated surface induces similar structures. In particular, 
given any hermitian inner product on the complex valued simplicial 
$1$-cochains,
the harmonic cochains split as holomorphic and anti-holomorphic $1$-cochains.
We'll prove analogues of the bilinear relations of Riemann, 
and show how one obtains a combinatorial period matrix. This construction
yields it's own combinatorial Schottky problem, 
which won't discuss in the current paper.

After describing our combinatorial construction, 
we'll show that if the complex valued simplicial cochains of a triangulated 
orientable Riemannian 2-manifold are equipped with a particular inner product 
induced by the Whitney embedding, then all of these structures 
provide a good approximation to the their continuum analogues.
In particular, the holomorphic and anti-holomorphic $1$-cochains
converge to the holomorphic and anti-holomorphic $1$-forms,
and the combinatorial period matrix converges to the conformal period
matrix of the associated Riemann surface, 
as the mesh of the triangulation tends to zero.
Hence, every conformal period matrix is a limit point of a sequence of 
combinatorial period matrices.

These statements may be interpreted as saying that a
triangulation of a surface, endowed with an inner product on the associated 
cochains, determines a conformal structure. Furthermore, for triangulations
of a Riemannian 2-manifold, a conformal structure is recovered (in the
limit) from algebraic and combinatorial data. 
Statements like this have been expressed by 
physicists for some time in various field theories and statistical mechanics.

We now proceed to describe the construction of combinatorial period matrices. 
First we need to extend some of our definitions from previous 
sections to deal with complex valued cochains. 
Let $\langle , \rangle$ be any non-degenerate positive definite hermitian 
inner product on the complex valued simplicial $1$-cochains
of a triangulated topological surface $K$. We define the associated 
combinatorial star operator $\bigstar$ by:
$$
\langle \bigstar a , b \rangle = (a \cup \overline{b})[M],
$$
where the bar denotes complex conjugation and $\cup$ is as in 
section~\ref{sec:prod}, extended linearly over $\mathbb{C}$. 
Just as with real coefficients, we have a Hodge decomposition
$$
C^1(K) \cong \delta C^0(K) \oplus H^1(K) \oplus \delta^{\ast} C^2(K)
$$
where $H^1$ is the space of complex valued harmonic $1$-cochains. 
Since $ \delta^{\ast}
\bigstar = \bigstar \delta$, by theorem ~\ref{thm:starprop}, 
$\bigstar$ induces an isomorphism of $H^1$. 

\begin{defn} Let $K$, $\langle , \rangle$, and $\bigstar$ be as above. 
We define subspace of the holomorphic $1$-cochains by
$$
\mathcal{H}^{1,0}(K) = \{ \sigma \in H^1(K) | \bigstar \sigma = -i \lambda 
\sigma \quad \textrm{for some} \quad \lambda \geq 0 \}
$$
and the subspace of the anti-holomorphic $1$-cochains by
$$
\mathcal{H}^{0,1}(K) = \{ \sigma \in H^1(K) | \bigstar \sigma = i \lambda 
\sigma \quad \textrm{for some} \quad \lambda \geq 0 \}
$$
\end{defn}

Since $\bigstar$ is not an orthogonal map, $\lambda$ may not equal one.
The following theorem shows that the space of harmonic $1$-cochains splits
into the subspaces of holomorphic and anti-holomorphic cochains.

\begin{thm} \label{thm:Hsplit}  Let $K$ be a triangulation of 
a surface $M$ of genus $g$ with
with a hermitian inner product on the simplicial 1-cochains of $K$, 
and the induced operator $\bigstar$.  
There is an orthogonal direct sum decomposition:
$$
H^1(K) \cong \mathcal{H}^{1,0} \oplus \mathcal{H}^{0,1}
$$
and each summand on the right has complex dimension $g$. Furthermore,
complex conjugation maps $\mathcal{H}^{1,0}$ to  $\mathcal{H}^{0,1}$ and
vice versa.
\end{thm}

\begin{proof}  The last assertion follows since $\bigstar$ is linear over 
$\mathbb{C}$. To prove the decomposition, we first note that the induced map 
$\bigstar$ on $\mathcal{H}$ has pure imaginary eigenvalues since
it is skew-adjoint: $\langle \bigstar \sigma , \tau \rangle =
- \langle \sigma , \bigstar \tau \rangle$.  
If $\sigma_1 \in \mathcal{H}^{1,0}$ and $\sigma_2 \in \mathcal{H}^{0,1}$ 
then for some 
$\lambda_1 , \lambda_2 > 0$ then
$$
-i \lambda_1 \langle \sigma_1 , \sigma_2 \rangle
= \langle \bigstar \sigma_1 , \sigma_2 \rangle
= - \langle \sigma_1 , \bigstar \sigma_2 \rangle
= i \lambda_2 \langle \sigma_1 , \sigma_2 \rangle
$$
so $\mathcal{H}^{1,0}$ and $\mathcal{H}^{0,1}$ are orthogonal. Finally,
$dim(\mathcal{H}^{1,0}) = dim(\mathcal{H}^{1,0}) = g$ since
$dim(\mathcal{H})= 2g$ and the eigenvalues of $\bigstar$ are 
all non-zero and occur in conjugate pairs.
\end{proof}

We'll now study further properties of holomorphic and anti-holomorphic
$1$-cochains. As in the smooth case, there is much to be gained by the 
analyzing the periods of these cochains. To do this, we first give
a brief description of the homology basis we'll evaluate these cochains on.

Without loss of generality, we assume that closed surface $M$ of genus $g$ is
obtained by identifying sides of a $4g$-gon, as in the following figure:
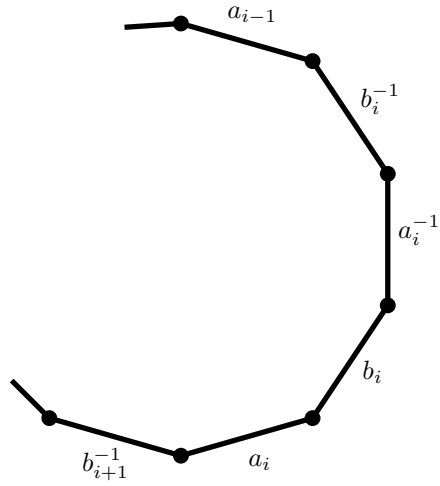
\begin{center}
\begin{pspicture}(0,.3)(8,7.4) 
\psline[linewidth=2pt,showpoints=true]{*-*}(1,1.5)(2.75,1)(4.5,1.5)(5.5,3)
(5.5,4.75)(4.5,6.25)(2.75,6.75)
\psline[linewidth=2pt]{-}(.5,2)(1,1.5)
\psline[linewidth=2pt]{-}(2.75,6.75)(2,6.7)
\psset{labelsep=4pt}
\uput[d](1.75,1.25){$b_{i+1}^{-1}$}
\uput[d](3.8,1.15){$a_i$}
\uput[d](5.3,2.4){$b_i$}
\uput[r](5.5,4){$a_{i}^{-1}$}
\uput[u](5.4,5.4){$b_{i}^{-1}$}
\uput[u](3.7,6.6){$a_{i-1}$}
\end{pspicture}
\end{center}
The basis $\{ a_1, a_2, \ldots a_g, b_1, b_2, \ldots b_g \}$ for the first 
homology is classically referred to as 
the canonical basis \cite{FK}, \cite{GS}, since it satisfies the following 
nice property: the intersection of any two basis elements 
is non-zero only for $a_i$ and $b_i$, in which case
it equals one. Of course, this basis is not canonical; nevertheless, we'll
work with it. (As a note, the discussion below is basis independent up to
an action of the modular symplectic group; we'll not go into details here.)
We assume our triangulation $K$ is a subdivision of the cellular 
decomposition given by the canonical homology basis. For any such
subdivision, each element of the canonical homology basis is represented 
as a sum of the edges into which it is subdivided, as in the following 
figure:
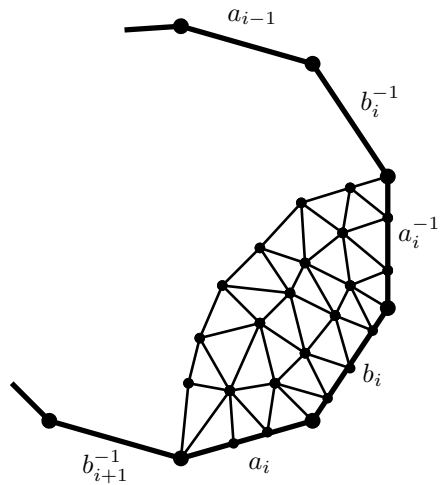
\begin{center}
\begin{pspicture}(0,.3)(8,7.4)
\psline[linewidth=2pt,showpoints=true]{*-*}(1,1.5)(2.75,1)(4.5,1.5)(5.5,3)
(5.5,4.75)(4.5,6.25)(2.75,6.75)
\psline[linewidth=2pt]{-}(.5,2)(1,1.5)
\psline[linewidth=2pt]{-}(2.75,6.75)(2,6.7)
\psset{labelsep=4pt}
\uput[d](1.75,1.25){$b_{i+1}^{-1}$}
\uput[d](3.8,1.15){$a_i$}
\uput[d](5.3,2.4){$b_i$}
\uput[r](5.5,4){$a_{i}^{-1}$}
\uput[u](5.4,5.4){$b_{i}^{-1}$}
\uput[u](3.7,6.6){$a_{i-1}$}
\psline[linewidth=1pt, showpoints=true]{*-*}(3.9,1.35)(4,2)(4.4,2.4)(4.8,2.9)
(5,3.3)(5.5,3.5)
\psline[linewidth=1pt, showpoints=true]{*-*}(4.5,1.5)(4,2)(4.7,1.8)
(4.4,2.4)(5,2.2)(4.8,2.9)(5.3,2.7)(5,3.3)(5.5,3)(5.5,3.5)
\psline[linewidth=1pt, showpoints=true]{*-*}(3.45,1.2)(3.4,1.9)(3.8,2.8)
(4.2,3.2)(4.4,3.6)(4.9,4)(5.5,4.2)
\psline[linewidth=1pt, showpoints=true]{*-*}(3.9,1.35)(3.4,1.9)(4,2)(3.8,2.8)
(4.4,2.4)(4.2,3.2)(4.8,2.9)(4.4,3.6)(5,3.3)(4.9,4)(5.5,3.5)
\psline[linewidth=1pt, showpoints=true]{*-*}(2.75,1)(2.85,2)(3,2.6)(3.3,3.3)
(3.8,3.8)(4.35,4.4)(5,4.6)(5.5,4.75)
\psline[linewidth=1pt, showpoints=true]{*-*}(2.75,1)(3.4,1.9)(2.85,2)
\psline[linewidth=1pt, showpoints=true]{*-*}(3.4,1.9)(3,2.6)(3.8,2.8)(3.3,3.3)
(4.2,3.2)(3.8,3.8)(4.4,3.6)(4.35,4.4)(4.9,4)(5,4.6)(5.5,4.2)
\end{pspicture}
\end{center}
By evaluating a cochain of $K$ on an element of the
canonical homology basis, we mean evaluating it on this subdivided 
representative.

\begin{defn} For $h \in \mathcal{H}$, the periods A-periods and B-periods
of $h$ are the following complex numbers:
$$
A_i = h(a_i) \quad  B_i = h(b_i) \qquad \textrm{for} \quad 1 \leq i \leq g
$$
\end{defn}

\begin{thm} \label{thm:rbr} [Riemann's Bi-linear relations] If $\sigma, 
\sigma' \in \mathcal{H}^{1,0}$, then the A-periods and B-periods satisfy:
$$
-i \lambda \langle \sigma , \overline{\sigma'} \rangle = \sum_{i=1}^g
(A_i B_i' - B_i A_i') = 0
$$
where $\lambda$ is such that $\bigstar \sigma = -i \lambda \sigma$.
\end{thm}

\begin{proof} Since $\sigma' \in \mathcal{H}^{1,0}$, 
$\overline{\sigma'} \in \mathcal{H}^{0,1}$ it follows that 
$\langle \sigma , \overline{\sigma'} \rangle = 0$. To show the
bi-linear relation we compute:
$$
-i \lambda \langle \sigma , \overline{\sigma'} \rangle = 
\langle \bigstar \sigma , \overline{\sigma'} \rangle =
(\sigma \cup \sigma')[M]
$$
where the fundamental class $[M]$ of $M$ may be represented by the sum
of the $2$-cells of $K$ appropriately oriented.
Now let $p:U \to M$ be the universal cover, with $U$ triangulated so that 
$p$ is locally a linear isomorphism onto the triangulation $K$ of $M$. Let $S$
denote a fundamental domain in the triangulation of $U$ so that the
induced map $p_{\ast}$ maps the $2$-simplicies of $S$ isomorphically onto
the $2$-simplicies of $K$. Then  $p_{\ast}(S) = [M]$, so the last 
expression equals
$$
(\sigma \cup \sigma')([M]) =
(p^{\ast} \sigma \cup p^{\ast} \sigma')(S)
$$
where $p^{\ast}$ denotes the pull back on cohomology. Since $\sigma$ is 
holomorphic, it is closed, as is $p^{\ast} \sigma$. Since $\overline{S}$ is
contractible to a point, the restriction of $p^{\ast} \sigma$ to 
$\overline{S}$ may be written as $p^{\ast} \sigma = \delta f$ for some
$0$-cochain $f$. Thus, since $\delta \sigma' = 0$ we have:
{\setlength\arraycolsep{2pt}
\begin{eqnarray*}
-i \lambda \langle \sigma , \overline{\sigma'} \rangle & = &
(\delta f \cup p^{\ast} \sigma')(S) \\
& = & (f \cup p^{\ast} \sigma')(\partial S) \\
& = & \sum_{i=1}^g (f \cup p^{\ast} \sigma') (a_i + a_i^{-1} + b_i + b_i^{-1})
\end{eqnarray*}}
It remains to show that this last expression equals $\sum_{i=1}^g 
(A_i B_i' - B_i A_i' )$. To do this, we first derive a simple
relation for the values of $f$ on the $0$-simplicies contained in the
cycles of the canonical homology basis. Consider the following 
figure:
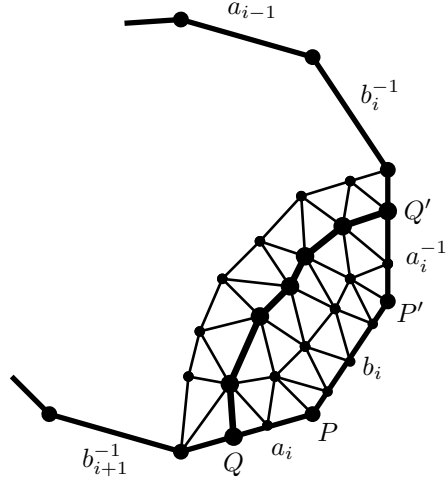
\begin{center}
\begin{pspicture}(0,.3)(8,7.4)
\psline[linewidth=2pt,showpoints=true]{*-*}(1,1.5)(2.75,1)(4.5,1.5)(5.5,3)
(5.5,4.75)(4.5,6.25)(2.75,6.75)
\psline[linewidth=2pt]{-}(.5,2)(1,1.5)
\psline[linewidth=2pt]{-}(2.75,6.75)(2,6.7)
\psset{labelsep=4pt}
\uput[d](1.75,1.25){$b_{i+1}^{-1}$}
\uput[d](4.1,1.25){$a_i$}
\uput[d](5.3,2.4){$b_i$}
\uput[r](5.6,3.6){$a_{i}^{-1}$}
\uput[u](5.4,5.4){$b_{i}^{-1}$}
\uput[u](3.7,6.6){$a_{i-1}$}
\psset{labelsep=6pt}
\uput[d](3.45,1.2){$Q$}
\uput[d](4.7,1.6){$P$}
\uput[d](5.8,3.2){$P'$}
\uput[r](5.5,4.2){$Q'$}
\psline[linewidth=1pt, showpoints=true]{*-*}(3.9,1.35)(4,2)(4.4,2.4)(4.8,2.9)
(5,3.3)(5.5,3.5)
\psline[linewidth=1pt, showpoints=true]{*-*}(4.5,1.5)(4,2)(4.7,1.8)
(4.4,2.4)(5,2.2)(4.8,2.9)(5.3,2.7)(5,3.3)(5.5,3)(5.5,3.5)
\psline[linewidth=2.5pt, showpoints=true]{*-*}(3.45,1.2)(3.4,1.9)(3.8,2.8)
(4.2,3.2)(4.4,3.6)(4.9,4)(5.5,4.2)
\psline[linewidth=1pt, showpoints=true]{*-*}(3.9,1.35)(3.4,1.9)(4,2)(3.8,2.8)
(4.4,2.4)(4.2,3.2)(4.8,2.9)(4.4,3.6)(5,3.3)(4.9,4)(5.5,3.5)
\psline[linewidth=1pt, showpoints=true]{*-*}(2.75,1)(2.85,2)(3,2.6)(3.3,3.3)
(3.8,3.8)(4.35,4.4)(5,4.6)(5.5,4.75)
\psline[linewidth=1pt, showpoints=true]{*-*}(2.75,1)(3.4,1.9)(2.85,2)
\psline[linewidth=1pt, showpoints=true]{*-*}(3.4,1.9)(3,2.6)(3.8,2.8)(3.3,3.3)
(4.2,3.2)(3.8,3.8)(4.4,3.6)(4.35,4.4)(4.9,4)(5,4.6)(5.5,4.2)
\end{pspicture}
\end{center}
The chain $\alpha$ from $Q$ to $Q'$ is a cycle. Since is homologous to the
cycle made up of chains from $Q$ to $P$, $P$ to $P'$ and $P'$ to $Q'$, and 
since the first and third push forward to the same chains on $K$, we have
that
$$
f(Q) - f(Q') = f( \partial \alpha) = \delta f ( \alpha ) = 
p^{\ast} \sigma (\alpha) = p^{\ast} \sigma (b_i) = B_i
$$
which means that for any $1$-cochains $p^{\ast} \tau$
$$
(f \cup \tau) (a_i^{-1}) = -((f + B_i) \cup \tau)(a_i)
= -(f \cup \tau) (a_i) - B_i \tau (a_i)
$$
Similarly,
$$
(f \cup \tau) (b_i^{-1}) = -((f - A_i) \cup \tau)(b_i)
= -(f \cup \tau) (b_i) + A_i \tau (a_i)
$$
So, we finally have that
{\setlength\arraycolsep{2pt}
\begin{eqnarray*}
-i \lambda \langle \sigma , \overline{\sigma'} \rangle & = &
\sum_{i=1}^g (f \cup p^{\ast} \sigma') (a_i + a_i^{-1} + b_i + b_i^{-1}) \\
& = & \sum_{i=1}^g -B_i p^{\ast} \sigma' (a_i) + A_i p^{\ast} \sigma' (b_i) \\
& = & \sum_{i=1}^g (A_i B_i' - B_i A_i' )
\end{eqnarray*}}
\end{proof}

Replacing $\overline{\sigma'}$ with $\sigma'$ in the previous proof shows
if $\sigma, \sigma' \in \mathcal{H}^{1,0}$ then
$$
-i \lambda \langle \sigma , \sigma' \rangle = \sum_{i=1}^g
(A_i \overline{B_i'} - B_i \overline{A_i'})
$$
where $\bigstar \sigma = -i \lambda \sigma$. If we apply this to 
$\sigma' = \sigma$ we obtain an expression for the norm of a holomorphic 
$1$-cochain in terms of its periods.

\begin{cor} \label{cor:rbrc} If $\sigma, \in \mathcal{H}^{1,0}$ satisfies  
$\bigstar \sigma = -i \lambda \sigma$ with periods $A_i$ and $B_i$ then
$$
\lVert \sigma \rVert^2 = \langle \sigma , \sigma \rangle =
\frac{i}{\lambda} \sum_{i=1}^g
(A_i \overline{B_i} - B_i \overline{A_i}) \geq 0
$$
\end{cor}

From this corollary we immediately have:

\begin{cor} Let $\sigma$ be a holomorphic $1$-cochain.
\begin{enumerate}
\item If the $A$-periods or $B$-periods of $\sigma$ vanish then $\sigma = 0$.
\item If the $A$-periods and $B$-periods of $\sigma$ are real then 
$\sigma = 0$.
\end{enumerate}
\end{cor}

Now let $\{ \tau_1, \tau_2, \ldots , \tau_g \}$ be a basis for the space of 
holomorphic cochains. By the corollary, if all the $A$-periods of a
linear combination of this basis vanish, then this linear combination 
is identically zero.
This implies we can solve uniquely for coefficients $c_{i,j}$ such that:
$$
\sum_{i=1}^g c_{i,j} \ \tau_i(a_k) = \delta_{j,k} 
$$

We put $\sigma_j = \sum_{i=1}^g c_{i,j} \ \tau_i$ and we call the basis 
$ \{ \sigma_1, \sigma_2, \ldots , \sigma_g \} $ \emph{the canonical basis 
for the space of holomorphic $1$-cochains}. From this we obtain the following
array of periods:

$$
\begin{array}{c|cccccccc}
& a_1 & a_2 & \cdots & a_g & b_1 & b_2 & \cdots & b_g \\
\hline
\sigma_1 & 1 & 0 & \cdots & 0 & \sigma_1(b_1) & \sigma_1(b_2)  
& \cdots & \sigma_1(b_g) \\
\sigma_2 & 0 & 1 & \cdots & 0 & \sigma_2(b_1) & \sigma_2(b_2)  
& \cdots & \sigma_2(b_g) \\
\vdots & & & \ddots & & & & & \\
\sigma_g & 0 & 0 & \cdots & 1 & \sigma_g(b_1) & \sigma_g(b_2)  
& \cdots & \sigma_g(b_g) \\
\end{array}
$$

\begin{defn} Let $ \{ \sigma_1, \sigma_2, \ldots , \sigma_g \} $ be 
the canonical basis for the space of holomorphic $1$-cochains and 
 $\{ a_1, a_2, \ldots a_g, b_1, b_2, \ldots b_g \}$ the canonical homology
basis, so $\sigma_i (a_j) = \delta_{i,j}$. We define the period 
matrix $\Pi = ( \pi_{i.j} )$ to be the $g \times g$ matrix of $B$-periods:
$$
\pi_{i,j} = \sigma_i (b_j)
$$
When we wish to emphasize the dependence of $\Pi$ on $K$ or $\langle,\rangle$
we'll write $\Pi_{K}$ or $\Pi_{K,\langle,\rangle}$. 
\end{defn}

\begin{rmk} Let $K$ be fixed. If two inner products on $C^1(K)$ differ by
a constant multiple then the associated period matrices are equal.
Hence, the combinatorial period matrix is a ``conformal invariant''.
\end{rmk}

\begin{thm} Let $K$ be a triangulated closed surface with a simplicial cochain
inner product. The associated period matrix $\Pi$ is symmetric and 
$Im( \Pi)$ is positive definite.
\end{thm}

\begin{proof} It suffices to show for $1 \leq i,j \leq g$, 
$\sigma_i(b_j) = \sigma_j(b_i)$ where $\sigma_i$ and $\sigma_j$ are 
canonical holomorphic cochain basis elements.
We apply theorem~\ref{thm:rbr} and compute:
{\setlength\arraycolsep{2pt}
\begin{eqnarray*}
0 & = & -i \lambda_i \langle \sigma_i , \overline{\sigma_j} \rangle \\
& = & \sum_{k=1}^g \sigma_i(a_k) \sigma_j(b_k) - \sigma_i(b_k) \sigma_j(a_k)\\
& = &  \sum_{k=1}^g \delta_{i,k} \sigma_j(b_k) - \delta_{j,k} \sigma_i(b_k)\\
& = & \sigma_j(b_i) - \sigma_i(b_j)
\end{eqnarray*}}
To prove the second statement, let $\sigma = \sum_{i=1}^g c_i \sigma_i$
be a nontrivial $\mathbb{R}$-linear combination of the canonical basis 
of holomorphic cochains. Then $\sigma(a_i) = c_i$. We show 
$$
\sigma \cdot Im(\Pi) \cdot \sigma > 0
$$ 
by using corollary~\ref{cor:rbrc} and computing:
{\setlength\arraycolsep{2pt}
\begin{eqnarray*}
0 & < & \frac{i}{\lambda} \sum_{k=1}^g \sigma(a_k) \overline{\sigma(b_k)} 
- \sigma(b_k) \overline{\sigma(a_k)} \\
& = & \frac{i}{\lambda} \sum_{k=1}^g c_k \Big( \sum_{i=1}^g c_i \overline{
\sigma_i(b_k)} \Big)
- c_k\Big( \sum_{i=1}^g c_i \sigma_i(b_k) \Big) \\
& = & \frac{i}{\lambda} \sum_{i=1}^g \sum_{k=1}^g c_k c_i \overline{
\sigma_i(b_k)} - c_k c_i \sigma_i(b_k) \\
& = & \frac{2}{\lambda} \sum_{i=1}^g \sum_{k=1}^g c_k c_i \ Im (\sigma_i(b_k)) \\
& = & \frac{2}{\lambda} \ ( \sigma \cdot Im(\Pi) \cdot \sigma )
\end{eqnarray*}}
\end{proof}

To this point, we have assumed $K$ is a triangulated closed topological 
surface and 
$\langle, \rangle$ is a non-degenerate inner product on the simplicial 
cochains of $K$. As remarked in the beginning of this section, the structures 
we have uncovered (splitting of harmonics,
bilinear relations, period matrix etc.) also appear for 
$1$-forms on a Riemann surface. In fact, all of the statements proven above
hold for forms as well \cite{GS}, except one should set $\lambda = 1$, 
since in this case the Hodge star operator $\star$ is an orthogonal 
transformation.

Now let $M$ be an orientable closed Riemannian 2-manifold. The Riemannian
metric induces an operator $\star$ which squares to $-Id$, and by identifying
tangent and cotangent space via the metric, this operator $\star$ gives
an almost complex structure. It is a theorem of Gauss that $M$ admits a
unique complex structure, i.e. Riemann surface structure, that is 
compatible with this almost complex structure. This theorem is, a priori,
non-trivial, and involves a transcendental construction of holomorphic
coordinates charts. By Torelli's theorem the resulting complex structure is 
determined uniquely by the period matrix of the associated Riemann surface $M$.

We can extend the usual $\mathcal{L}_2$ inner product on the vector space of 
real valued $1$-forms to a hermitian inner product on 
the space of complexified $1$-forms canonically, by declaring 
\begin{equation} \label{eq:ip}
\langle \omega_1 \otimes z_1 , \omega_2 \otimes z_2 \rangle =
z_1 \overline{z_2} \langle \omega_1 , \omega_2 \rangle
\end{equation}
for $\omega_i \otimes z_i \in T^* M \bigotimes \mathbb{C}$. Let 
$\| \ \ \| $ denote the induced norm.

Now, for any triangulation $K$, by embedding complex valued $1$-cochains
into $T^* M \bigotimes \mathbb{C}$ via the Whitney embedding,
we obtain the induced Whitney inner product on complex valued $1$-cochains.
For the remainder of this section, we work only with this Whitney inner 
product.
We remark here that while the approximation theorems from 
section~\ref{sec:whitney} and ~\ref{sec:star} (using
the Whitney inner product) involved real-valued forms and cochains, 
the proofs follow verbatim for complex coefficients as well.

First we prove the following:

\begin{lemma} Let $M$ be a Riemannian 2-manifold with triangulation $K$ of 
mesh $\eta$, and $\mathfrak{h}$ a be complex valued $1$-form on $M$, so 
$\star \mathfrak{h}  = -i \mathfrak{h}$. By the 
Hodge decomposition of cochains and 
theorem~\ref{thm:Hsplit} we may write
$$
R \mathfrak{h} = \delta g + h_1 + h_2 + \delta^{\ast} k
$$
uniquely for $h_1 \in \mathcal{H}^{1,0}$ and $h_2 \in \mathcal{H}^{0,1}$.
Then there exists a positive constant $C$, independent of $K$, such that 
$$
\lVert W h_1 - \mathfrak{h} \rVert \leq C \cdot \eta
$$ 
\end{lemma}

\begin{proof}
By theorems ~\ref{thm:happ} and ~\ref{thm:happstar}, there is a 
positive constant $C$, independent of $K$, such that
{\setlength\arraycolsep{2pt}
\begin{eqnarray*}
C \cdot \eta & \geq & \| W \bigstar (h_1 + h_2) - \star \mathfrak{h} \| 
+ \| \mathfrak{h} - W (h_1 + h_2) \| \\
& = &  \| W \bigstar (h_1 + h_2) - \star \mathfrak{h} \|
+ \| \star \mathfrak{h} + i W (h_1 + h_2) \| \\
& \geq &  \| W \bigstar h_1 + W \bigstar h_2 + iW(h_1 + h_2) \| \\
& = &  \| \bigstar h_1 + \bigstar h_2 + i(h_1 + h_2) \| \\
\end{eqnarray*}}
Since $h_1 \in \mathcal{H}^{1,0}$ and $h_2 \in \mathcal{H}^{0,1}$ we may
write $\bigstar h_1 = -i \lambda_1 h_1$ and  $\bigstar h_2 = i \lambda_2 h_2$
for some $\lambda_1,\lambda_2 > 0 $. Using the fact that 
$\mathcal{H}^{1,0} \perp \mathcal{H}^{0,1}$ we then have
{\setlength\arraycolsep{2pt}
\begin{eqnarray*}
C^2 \cdot \eta^2 & \geq & 
\| -i \lambda_1 h_1 + i \lambda_2 h_2 + i h_1 +  i h_2 \|^2 \\
& = & \| (1 - \lambda_1) h_1 + (1 + \lambda_2) h_2  \|^2 \\
& = & \langle (1 - \lambda_1) h_1,  (1 - \lambda_1) h_1 \rangle +
\langle (1 + \lambda_2) h_2, (1 + \lambda_2) h_2 \rangle \\
& = & |1 - \lambda_1|^2 \|h_1\|^2 + |1 + \lambda_2|^2 \|h_2\|^2
\end{eqnarray*}}
So, we conclude
$$
\| h_2 \| \leq \frac{C \cdot \eta}{|1 + \lambda_2|} \leq C \cdot \eta
$$
and finally, 
$$
\| W h_1 - \mathfrak{h} \|  \leq 
\| W(h_1 + h_2) - \mathfrak{h} \rVert + \| h_2 \| \leq 2C \cdot \eta
$$
\end{proof}

\begin{rmk} A closer examination of the proof shows that, for a fine
enough triangulation, $1 - \lambda_1$ is bounded by a constant times the 
mesh.
There is, of course, an analogous statement for anti-holomorphic
1-forms $\mathfrak{h}$ and the anti-holomorphic part of the cochain 
$R \mathfrak{h}$.
\end{rmk}

One can check that the hermitian inner product on $1$-forms of $M$, 
defined in \eqref{eq:ip}, agrees with the usual inner product on the
$1$-forms of the Riemann surface associated to $M$, given by
$$
\langle \omega , \eta \rangle = \int_M \omega \wedge \star \overline{\eta}.
$$ 
It is a peculiarity of working in the middle dimension (here 1) that this
inner product, and the Hodge star operator, depend only on the conformal
class of the Riemannian metric. This implies that the period matrix of the
Riemann surface associated to $M$ can be computed by using the inner
product in \eqref{eq:ip} in the following way: 
split off the harmonic $1$-forms and 
evaluate the appropriate basis of the $-i$ eigenspace of $\star$ on the
canonical homology basis. We remark here that this involves
a transcendental procedure in the Hodge decomposition of forms. The point
of the following theorem is that the period matrix, and therefore the
complex structure, is computable to any desired accuracy, from algebraic
and combinatorial data.

\begin{thm} Let $M$ be a closed orientable Riemannian 2-manifold
and let $\Pi$ be the period matrix of the Riemann surface associated to $M$.
Let $K_n$ be a sequence of triangulations of $M$ with mesh converging to zero. 
Then, for each $n$, the induced Whitney 
inner product on the simplicial 1-cochains of $K_n$ gives rise to a 
combinatorial period matrix $\Pi_{K_n}$, and
\begin{displaymath}
\lim\limits_{n \to \infty} \Pi_{K_n} = \Pi.
\end{displaymath}
\end{thm} 

\begin{proof}
Let $\mathfrak{h}_g, \cdots , \mathfrak{h}_g$ be the
canonical basis of holomorphic $1$-forms with periods
{\setlength\arraycolsep{2pt}
\begin{eqnarray*}
\mathfrak{h}_i (a_j) &=& \int_{a_j} \mathfrak{h}_i \ = \ \delta_{i,j} \\
\mathfrak{h}_i (b_j) &=& \int_{b_j} \mathfrak{h}_i \ = \ \pi_{i,j} 
\end{eqnarray*}}
for $1 \leq i,j \leq g$, and $\pi_{i,j}$ the $(i,j)$ entry of $\Pi$.

For each $n$, let $\varphi_1^n, \cdots , \varphi_g^n$ be a
basis for the holomorphic cochains on $K_n$. Then the periods are
{\setlength\arraycolsep{2pt}
\begin{eqnarray*}
\varphi_i^n (a_j) & = & \delta_{i,j} \\
\varphi_i^n (b_j) & = & \pi_{i,j}^n
\end{eqnarray*}}
for $1 \leq i,j \leq g$, and $\pi_{i,j}^n$ the $(i,j)$ entry of $\Pi_{k_n}$.
Our goal is to show, for all $1 \leq i,j \leq g$,
$$
\lim\limits_{n \to \infty} \varphi_i^n (b_j) = \mathfrak{h}_i (b_j).
$$

Let $R_n$ denote the integration map taking 1-forms to cochains on $K_n$.
We define $h_i^n$ holomorphic part of the cochain $R_n \mathfrak{h}_i$.
By the previous lemma, $h_i^n \to \mathfrak{h}_i$ as $n \to \infty$. Therefore, by evaluating on a cycle $a_j$, we see from the Hodge decomposition of these closed forms that
\begin{equation} \label{eq:limit}
\lim\limits_{n \to \infty} h_i^n (a_j) = \mathfrak{h}_i (a_j) = \delta_{i,j}
\end{equation}
For each $n$ and $1 \leq i \leq g$ we may write 
$$
h_i^n = \sum_{k=1}^g c_{i,k}^n \varphi_k^n
$$
and by evaluating on the cycle $a_j$ we see similarly that
$$
c_{i,j}^n  = \sum_{k=1}^g c_{i,k}^n \varphi_k^n (a_j) = h_i^n (a_j)
$$
Combining this with equation \eqref{eq:limit}, we have
$$
\lim\limits_{n \to \infty} c_{i,j}^n = \delta_{i,j}
$$
which implies
$$
\lim\limits_{n \to \infty} \| \varphi_i^n - h_i^n \| = 0
$$
By the lemma, $\| h_i^n \| \to \| \mathfrak{h}_i \|$, so 
the sequences $\| h_i^n \|$ and $\| \varphi_i^n \|$ are bounded.  
Finally, we have
$$
\lim\limits_{n \to \infty} \varphi_i^n (b_j) = \lim\limits_{n \to \infty} 
h_i^n (b_j) = \mathfrak{h}_i (b_j)
$$
\end{proof}

\begin{cor}Let $M$ be a closed Riemann surface with period matrix $\Pi$.
Let $K_n$ be a sequence of triangulations of $M$ with mesh converging to zero,
and combinatorial period matrices $\Pi_{K_n}$ induced by the Whitney metric. 
Then, 
$$
\lim\limits_{n \to \infty} \Pi_{K_n} = \Pi.
$$
\end{cor}

\begin{proof} While there isn't a notion of geodesic length on a 
Riemann surface, a distance converging to zero is well defined since it
depends only on a conformal class of metrics. So the statement of the 
corollary makes sense. Then one can choose any Riemannian metric on $M$ in 
the conformal class of metrics determined by $M$, and apply the above theorem.
\end{proof}

\begin{cor} Every conformal period matrix is the limit of a sequence of 
combinatorial period matrices.
\end{cor}

\section{Inner Products and Their Inverses} \label{sec:ip}

In this section we study inner products on cochains, as well as the 
induced ``inverse inner product''. Smits also studied the inverse of 
inner products in \cite{}, where her proved his results on the convergence of
the divergence operator $d^{\ast}$ on a surface.

We start with the following definition:
\begin{defn} A \emph{geometric inner product} on the simplicial
cochains $C = \bigoplus_j C^j$ of a triangulated space $K$ 
is a non-degenerate positive definite inner product $\langle , \rangle$ on $C$ 
satisfying:
\begin{enumerate}
\item $C^i \perp C^j$ for $i \neq j$
\item locality: $\langle a , b \rangle \neq 0$ only if $St(a) \cap St(b)$
is non-empty.
\end{enumerate}
\end{defn}

\begin{rmk} A geometric inner product restricted to 1-cochains, and its
induced norm, gives a notion of lengths of edges and the angles between them.
It may be interesting to study the consequences of an inner product of 
signature other than the one considered here.
\end{rmk}

We assume in this section that all cochain inner products are geometric in the
above sense. Note that the Whitney inner product is geometric.

An inner product on $C^*$ induces an isomorphism from $C^*$
to the linear dual of $C^*$, which we denote by $C_*$ and refer to as the 
simplicial chains (to be more precise, this is the double dual of chains, 
but we'll confuse the two since we're assuming $K$ is compact). 
The inverse of the inner product is, by definition,
the inverse of the isomorphism $C^* \to C_*$, and is an isomorphism
$C_* \to C^*$. This gives an inner product on the (simplicial) chains $C_*$ 
and will be denoted by $\langle , \rangle^{-1}$. 

If one represents a geometric cochain inner product as a matrix, using
the standard basis given by the simplicies, then the 
locality property roughly states that this matrix is ``near diagonal''. 
Of course, the inverse of a diagonal matrix is diagonal, but the inverse of a 
near diagonal matrix is \emph{not} near diagonal. Rather, it can have all
entries non-zero; i.e. the inverse inner product on chains is \emph{not}
geometric.\footnote{It is true is that the matrix entries decrease in 
absolute value as they move from the diagonal, so that the inner product 
of two chains decays rapidly as a function of ``geometric distance''.}

In this section, we describe the inner product  $\langle , \rangle^{-1}$
on chains in a geometric way
by showing it can be expressed as a weighted sum of paths 
in a collection of graphs associated to $K$. This will be useful in the 
next section for making explicit computations of the combinatorial star
operator. We begin with some definitions:

\begin{defn}A graph $\Gamma$ (without loops) consists of a set $S$, 
called vertices, and a 
collection of cardinality two subsets of $S$, called edges.
Two edges of $\Gamma$ are said to be incident if their intersection (as 
subsets of $S$) is nonempty.
A weighted graph is a graph with an assignment of a real number $w(e)$ 
to each edge $e$.
A path $\gamma$ in a graph is a sequences of edges $\{e_i\}_{i \in I}$ such 
that for each $i$, $e_i$ and $e_{i+1}$ are incident to a common vertex.
The weight $w(\gamma)$ of a path $\gamma$ in a weighted graph 
is the product of the weight of each edge in $\gamma$.
By convention, we say there is a unique path of length zero between
any vertex and itself, and the weight of this path is one.
\end{defn}

\begin{defn}Let $K$ be the simplicial cochain complex of a triangulated 
n-manifold $M$. We define the  
\emph{graph associated to the $j$-simplicies of $K$}, denoted $\Gamma(K,j)$,
to be the following graph: The vertices of $\Gamma(K,j)$ are the 
set $\{\sigma_{\alpha}\}$ of $j$-simplicies of $K$; two distinct vertices 
$\sigma_1, \sigma_2$ of $\Gamma(K,j)$ form an edge 
if and only if they are faces of a common $n$-simplex of $K$ (i.e 
$St(\sigma_1) \cap St( \sigma_2)$ is non-empty).
\end{defn}

\begin{cor} \label{cor:graph}Let $K$ be the simplicial cochain complex of a 
triangulated n-manifold $M$.
\begin{enumerate}
\item
Paths in $\Gamma(K,j)$ correspond to sequences $\{s_i\}_{i \in I}$
of $j$-simplicies in $K$ such that, for each $i$, $s_i$ and $s_{i+1}$ are
faces of a common $n$-simplex.
\item $\Gamma(K,0)$ is isomorphic to $K_1$, the $1$-skeleton of $K$ (the union
of its vertices and edges).
\end{enumerate}
\end{cor}

\begin{proof}
This follows since $K$ is homeomorphic to a manifold.
\end{proof}

Now suppose the cochains $C^*$ of $K$ are endowed with a geometric 
inner product $\langle, \rangle$. (Our
motivating example is the Whitney metric on $C^*$, but other examples 
arise when considering interactions on simplicial lattices.) 
In this case we associate to $(C^* , \langle ,\rangle)$ the following 
collection of weighted graphs.

\begin{defn}Let $C^*$ be the cochains of a finite triangulation $K$ of a 
manifold, with geometric cochain inner product $\langle, \rangle$. 
We define the  
\emph{weighted graph associated to the $j$-cochains $C^j$}, denoted 
$\Gamma_w(K,j)$,
to be the following weighted graph: The underlying graph of 
$\Gamma_w(K,j)$ is $\Gamma(K,j)$ and the weight $w(e)$ of an edge 
$e = \{\sigma_1,\sigma_2\}$ equals 
$$
w(e) = \frac{\langle \sigma_1 , \sigma_2 \rangle}{\| \sigma_1 \| \cdot
\| \sigma_2 \|}
$$
where $\| \sigma \| = \sqrt{ \langle \sigma, \sigma \rangle}$
\end{defn}

\begin{rmk}  \label{rmk:graphrmk}
The appropriate analogue corollary ~\ref{cor:graph} for 
weighted graphs holds as well.
\end{rmk}

The following describes how the metric $\langle , \rangle^{-1}$ on 
$C_j$ can be computed by counting weighted paths in the 
weighted graph $\Gamma_w(K,j)$ associated to $(C^j \langle,\rangle)$.

\begin{thm}For $\sigma_1 , \sigma_2 \in C_j$
$$
\langle \sigma_1 , \sigma_2 \rangle ^{-1} =
\frac{1}{\|\sigma_1 \| \cdot \| \sigma_2 \|} \sum_{i \geq 0}
(-1)^i \sum_{\gamma_i \in \Gamma_w(K,j)} w(\gamma_i)
$$
where $\gamma_i$ is a path in $\Gamma_w(K,j)$ of length $i$, starting at
$\sigma_1$ and ending at $\sigma_2$.
\end{thm}

\begin{proof}
Let $M$ be the matrix for $\langle , \rangle$ with respect to a fixed
ordering of the basis given by the simplicies of $K$. Let $D$ be 
the diagonal matrix, with respect to the same ordered basis, 
whose diagonal entries are the norm of a simplex.
Let $|M| = D^{-1} M D^{-1}$. Note that the entries of $|M|$ are 
normalized since the entries of $D^{-1}$ are of the form 
$\frac{1}{\| \sigma \|}$. In particular the diagonal entries of 
$|M|$ equal $1$, so we may write 
$$
M^{-1} = D^{-1} |M|^{-1} D^{-1} =  D^{-1} (I + A)^{-1} D^{-1}
$$
It is easy to check that $A$ is precisely the weighted adjacency matrix 
for the weighted
graph $\Gamma_w(K,j)$. Recall that the $i^{th}$ power of a weighted adjacency 
matrix counts the sum of the weights of all paths of length $i$. 
Then, by the Cauchy-Schwartz inequality, the formula
$$
(I + A)^{-1} = \sum_{i \geq 0} (-1)^i A^i
$$
may be applied above, and we conclude that
$$
\langle \sigma_1 , \sigma_2 \rangle ^{-1} = \sigma_1 M^{-1} \sigma_2 = 
\frac{1}{\|\sigma_1 \| \cdot \| \sigma_2 \|} \sum_{i \geq 0}
(-1)^i \sum_{\gamma_i \in \Gamma_w(K,j)} w(\gamma_i)
$$
\end{proof}

\begin{rmk}
\begin{enumerate}
\item The above theorem in the case $j=0$, in light of 
remark ~\ref{rmk:graphrmk}, shows that for vertices $p$ and $q$ of $K$,
$\langle p , q \rangle^{-1}$ may be expressed
as a weighted sum over all paths in the $1$-skeleton $K_1 \subset K$.
\item These expressions for $\langle , \rangle ^{-1}$ not only provide
a nice geometric interpretation, but are also useful for computations,
as we will see in section~\ref{sec:circle} 
where we compute $\bigstar$ for the circle.
\end{enumerate}
\end{rmk}

\section{Computation for $S^1$} \label{sec:circle}
In this section we compute the operator $\bigstar$ explicitly for the circle
$S^1$. We take $S^1$ to the be the unit interval $[0,1]$ with $0$ and 
$1$ identified. We consider a
sequence of subdivisions, the $n^{\mathrm{th}}$ triangulation being given by
vertices at the points $v_i = \frac{i}{n}$ for $0 \leq i \leq n$. We
denote the edge from $v_i$ to $v_{i+1}$ by $e_i$ for $0 \leq i \le n$ and 
orient this edge from $v_i$ to $v_{i+1}$. See the following figure:
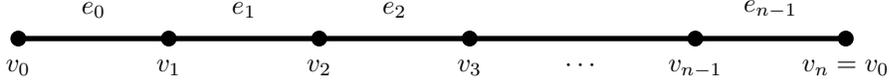
\begin{center}
\begin{pspicture}(0,0)(11,2)
\psline[linewidth=2pt,showpoints=true]{*-*}(0,1)(2,1)(4,1)(6,1)(9,1)(11,1)
\psset{labelsep=8pt}
\uput[u](1,1){$e_0$}
\uput[u](3,1){$e_1$}
\uput[u](5,1){$e_2$}
\uput[u](10,1){$e_{n-1}$}
\uput[d](0,1){$v_0$}
\uput[d](2,1){$v_1$}
\uput[d](4,1){$v_2$}
\uput[d](6,1){$v_3$}
\uput[d](7.5,1){$\cdots$}
\uput[d](9,1){$v_{n-1}$}
\uput[d](11,1){$v_n = v_0$}
\end{pspicture}
\end{center}

All operators will be written as matrices with respect to the ordered basis
$\{ v_0,\ldots,v_{n-1},e_0,\ldots,e_{n-1} \}$. 

Recall that the operator $\bigstar$ is defined by $\langle \bigstar \sigma, 
\tau \rangle  = (\sigma \cup \tau)[S^1]$ where here $[S^1]$ is the sum
of all the edges with their chosen orientations. 
We'll use the cochain inner product $\langle , \rangle$ induced by the Whitney
embedding and the standard metric on $S^1$ (i.e.  $\langle \mathrm{d}t, 
\mathrm{d}t \rangle = 1$.
Let $M$ denote the matrix for the cochain inner product and let $C$ denote 
the matrix for the pairing given by 
$(\sigma,\tau) \mapsto (\sigma \cup \tau)[S^1]$.
Then $\bigstar = M^{-1}C$. (We suppress the dependence of these 
operators on the level of subdivision; the $n^\mathrm{th}$ 
level $M$ and $C$ are size $2n \times 2n$.)

By the definition of $\cup$ and our chosen orientations we have that
\[
C =
\left( \begin{array}{c|c}
0 & A \\
\hline
A^t & 0
\end{array}
\right)
\]
where 
\[
A =
\left( 
\begin{array}
{cccccc}
1/2 & 0 & 0    & \ldots & 0 & 1/2 \\ 
1/2 & 1/2 & 0 & \ldots & \dots & 0 \\
0 & 1/2 & 1/2 & \ddots &       & \vdots \\
\vdots & \ddots & \ddots & \ddots & \ddots & \vdots \\
\vdots & & \ddots & \ddots & \ddots & 0 \\
0 & \ldots & \ldots & 0 & 1/2 & 1/2
\end{array}
\right)
\]
and $t$ denotes transpose.

One can compute explicitly:
\[
\langle \sigma, \tau \rangle = \left\{ \begin{array}{ll}
\frac{2}{3n} & \textrm{$\sigma = \tau$ is a vertex} \\
\frac{1}{6n} & \textrm{$\sigma$, $\tau$ are vertices in the boundary of a 
common edge} \\
n & \textrm{$\sigma = \tau$ is an edge} \\
0 & \textrm{otherwise}
\end{array} \right.
\]
So, in our chosen basis, the matrix for the inner product is given by:
\[
M =
\left( \begin{array}{c|c}
B & 0 \\
\hline
0 & nI
\end{array}
\right)
\]
where $I$ denotes the $n \times n$ identity matrix and 
\[
B =
\left( 
\begin{array}
{cccccc}
2/3n & 1/6n & 0    & \ldots & 0 & 1/6n \\ 
1/6n & 2/3n & 1/6n & \ldots & \dots & 0 \\
0 & 1/6n & 2/3n & \ddots &       & \vdots \\
\vdots & \ddots & \ddots & \ddots & \ddots & \vdots \\
0 & & \ddots & \ddots & \ddots & 1/6n \\
1/6n & 0 & \ldots & 0 & 1/6n & 2/3n
\end{array}
\right).
\]

We now compute $B^{-1}$. Note that one can 
write $B = \frac{2}{3n}(\frac{1}{4} D + I)$ where 
\[
D =
\left( 
\begin{array}
{cccccc}
0 & 1 & 0    & \ldots & 0 & 1 \\ 
1 & 0 & 1 & \ldots & \dots & 0 \\
0 & 1 & 0 & \ddots &       & \vdots \\
\vdots & \ddots & \ddots & \ddots & \ddots & \vdots \\
0 & & \ddots & \ddots & \ddots & 1 \\
1 & 0 & \ldots & 0 & 1 & 0
\end{array}
\right)
\]
then 
\begin{eqnarray*}
B^{-1} & = & \frac{3n}{2} \bigg( \frac{1}{4}D + I \bigg) ^{-1} \\
& = & \frac{3n}{2} ( I - \frac{1}{4}D + \frac{1}{4^2}D^2 - \frac{1}{4^3}D^3 
\pm \cdots) \\
& = & \frac{3n}{2} \sum_{k\geq 0 } (-1/4)^k D^k
\end{eqnarray*}

Note that $D$ is the adjacency matrix for the graph corresponding to the 
original triangulation $K$, or rather, $\frac{1}{4} D$ is the weighted 
adjacency matrix for the weighted graph in the following figure:

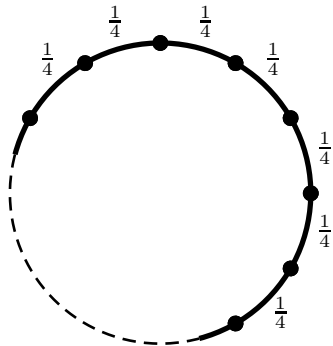
\begin{center}
\begin{pspicture}(0,-.3)(4.5,4.5)
\psarc[linewidth=2pt]{*-*}(2,2){2}{0}{30}
\psarc[linewidth=2pt]{*-*}(2,2){2}{30}{60}
\psarc[linewidth=2pt]{*-*}(2,2){2}{60}{90}
\psarc[linewidth=2pt]{*-*}(2,2){2}{90}{120}
\psarc[linewidth=2pt]{*-*}(2,2){2}{120}{150}
\psarc[linewidth=2pt]{*-}(2,2){2}{150}{165}
\psarc[linestyle=dashed, linewidth=1pt](2,2){2}{165}{285}
\psarc[linewidth=2pt]{-*}(2,2){2}{285}{300}
\psarc[linewidth=2pt]{*-*}(2,2){2}{300}{330}
\psarc[linewidth=2pt]{*-*}(2,2){2}{330}{360}
\uput[r](3.9,1.5){$\frac{1}{4}$}
\uput[r](3.9,2.6){$\frac{1}{4}$}
\uput[u](3.5,3.4){$\frac{1}{4}$}
\uput[u](2.6,3.9){$\frac{1}{4}$}
\uput[u](1.4,3.9){$\frac{1}{4}$}
\uput[d](3.6,.8){$\frac{1}{4}$}
\uput[u](.5,3.4){$\frac{1}{4}$}
\end{pspicture}
\end{center}

As shown in section~\ref{sec:ip}, the matrices $\frac{1}{4^k}D^k$ have a 
geometric interpretation: 
the $(i,j)$ entry equals the total weight of all paths  
from $v_i$ to $v_j$ of length $k$. Since in this case all weights are 
$\frac{1}{4}$, we'll simply compute the the $(i,j)$ entry of $D^k$, 
i.e the total number of paths from $v_i$ to $v_j$ of length $k$.

We first note that for the real line with integer vertices, the number of 
paths of length $r$ between two vertices distance $s$ apart is the 
binomial coefficient $\binom{r}{\frac{r+s}{2}}$. 
By considering the standard covering of the circle with $n$ vertices
by the line we have
\[
d_{i,j}^k = \sum_{t \in \mathbb{Z}} \binom{k}{\frac{k + \lvert i - j 
\rvert + nt}{2}}
\]
where the above binomial coefficient is zero unless 
$\frac{k + \lvert i - j \rvert + nt}{2}$ is a non-negative integer less
than or equal to $k$. Hence,
\[
M^{-1}=
\left( \begin{array}{c|c}
\bigg[  \frac{3n}{2} \displaystyle{\sum_{k \geq 0}} \bigg( \frac{-1}{4}\bigg)^k 
d_{i,j}^k \bigg] & 0 \\
\hline
0 & \frac{1}{n}I
\end{array}
\right).
\]
We conclude that:
\begin{eqnarray*}
\bigstar v_i & = & \frac{1}{2n}(e_{i-1} + e_i)  \\
\bigstar e_i & = & \frac{3n}{4} \sum_{0 \leq j \leq n-1} 
\bigg( \sum_{k \geq 0} \bigg( \frac{-1}{4} \bigg)^k 
\sum_{t \in \mathbb{Z}}  
\binom{k}{\frac{k + \lvert i - j \rvert + nt}{2}} 
+ \binom{k}{\frac{k + \lvert i - (j + 1) \rvert + nt}{2}}\bigg) v_j
\end{eqnarray*}

In the figures below, we plot $\bigstar e_{n/2}$ for $n=10,20,50$.
In each figure, the x-axis denotes the circle, triangulated with black
dots as vertices. For fixed $n$, and each $0 \leq i \leq n$, 
we plot the coefficient of $v_i$ appearing in $\bigstar e_{n/2}$. We've
used a triangle to denote this value. To suggest that the plots are roughly
Gaussian and that operator approaches a ``delta function'', 
we have connected consecutive plot points with a line.

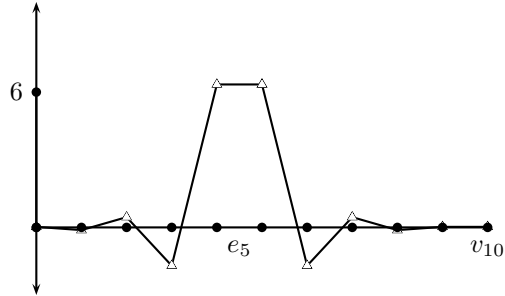
\begin{center}
\psset{xunit=.6cm,yunit=.3cm}
\begin{pspicture}(0,-5)(12,12) 
\savedata{\mydata}[
{{0, 0.02392344}, {1, -0.1196172}, {2, 0.4545455}, {3, -1.698565}, 
 {4, 6.339713},{5,6.339713}{5, 6.339713}, {6, -1.698565}, {7, 0.4545455}, 
{8, -0.1196172}, {9, 0.02392344}, {10, 0.02392344}} 
]
\dataplot[showpoints=true,dotstyle=triangle]{\mydata}
\savedata{\verts}[
{{0, 0}, {1, 0}, {2, 0}, {3, 0}, {4, 0}, {5, 0}, {6, 0}, {7, 0}, {8, 0}, {9, 
    0}, {10, 0}}
]
\dataplot[showpoints=true]{\verts}
\psline{<->}(0,10)(0,-3)
\psline{*-*}(0,0)(0,6)
\uput[d](4.5,0){$e_5$}
\uput[l](0,6){$6$}
\uput[d](10,0){$v_{10}$}
\end{pspicture}
\end{center}

\begin{center}
\psset{xunit=.4cm,yunit=.2cm}
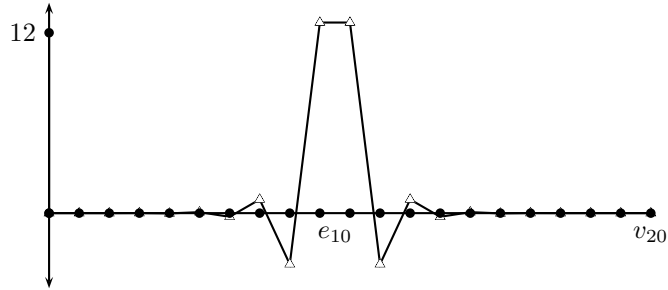
\begin{pspicture}(0,-5)(20,16) 
\savedata{\mydatatwent}[
{
{0, -0.00006608686}, {1,0.0003304343},{2, -0.001255650},
{3, 0.004692167},{4, -0.01751302},{5, 0.06535991},{6, -0.2439266},
{7,0.9103466},{8, -3.397460},{9, 12.67949},
{10, 12.67949}, {11, -3.397460}, {12,0.9103466}, {13, -0.2439266}, 
{14, 0.06535991}, {15, -0.01751302}, {16, 0.004692167}, 
{17, -0.001255650}, {18,0.0003304343}, {19, -0.00006608686},
{20, -0.00006608686}}
]
\dataplot[showpoints=true,dotstyle=triangle]{\mydatatwent}
\savedata{\vertstwent}[
{{0, 0}, {1, 0}, {2, 0}, {3, 0}, {4, 0}, {5, 0}, {6, 0}, {7, 0}, {8, 0}, {9, 
    0}, {10, 0}, {11, 0}, {12, 0}, {13, 0}, {14, 0}, {15, 0}, {16, 0}, {17, 
    0}, {18, 0}, {19, 0}, {20, 0}}
]
\dataplot[showpoints=true]{\vertstwent}
\psline{<->}(0,14)(0,-5)
\psline{*-*}(0,0)(0,12)
\uput[d](9.5,0){$e_{10}$}
\uput[d](20,0){$v_{20}$}
\uput[l](0,12){$12$}
\end{pspicture}
\end{center}

\begin{center}
\psset{xunit=.25cm,yunit=.12cm}
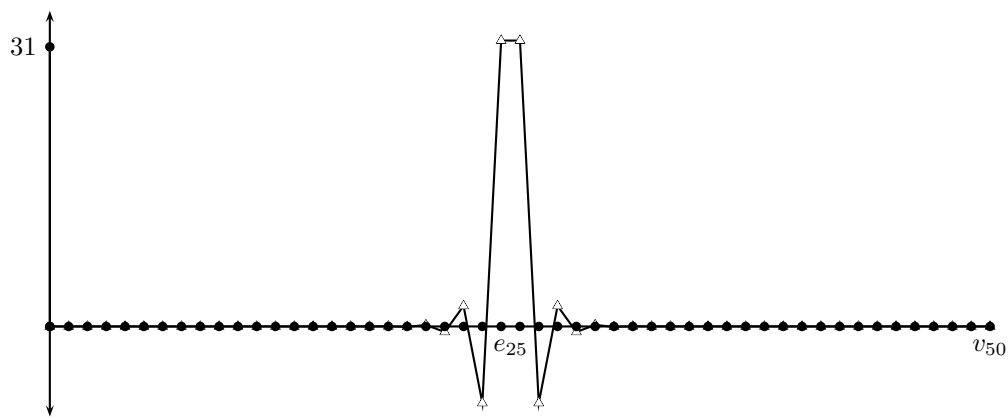
\begin{pspicture}(-2,-12)(50,40) 
\savedata{\moredata}[
{
{0,0.0000000000004353},
{1,-0.00000000000217},
{2,.00000000000827},
{3,-0.0000000000309},
{4,0.00000000011536},
{5,-0.0000000004305}
{6,.000000001606},
{7,-0.00000000599}, 
{8,.000000022381},
{9,-.000000083527},
{10,.00000031172},
{11,.00000116338}
{12,.0000043418},
{13,-0.0000162},
{14,0.0000604},
{15,-0.00022569},
{16,0.0008422},
{17,-0.00314346},
{18,0.0117315},
{19,-0.0437828},
{20,0.16339},
{21,-0.609816}, 
{22,2.2758664},
{23,-8.4936490},
{24, 31.698729},
{25, 31.698729}, 
{26,-8.4936490}, 
{27,2.2758664}, 
{28,-0.609816}, 
{29,0.16339},
{30,-0.0437828},
{31,0.0117315},
{32,-0.00314346}, 
{33,0.0008422},{34,-0.00022569},{35,0.0000604},{36,-0.0000162},
{37,.0000043418},{38,.00000116338},{39,.00000031172},
{40,-.000000083527},{41,.000000022381},{42,-0.00000000599}, 
{43,.000000001606},{44,-0.0000000004305},{45,0.00000000011536}, 
{46,-0.0000000000309},{47,.00000000000827},{48,-0.00000000000217}, 
{49,0.0000000000004353},{50,0.0000000000004353}
}
]
\dataplot[showpoints=true,dotstyle=triangle]{\moredata}
\savedata{\moreverts}[
{{0, 0}, {1, 0}, {2, 0}, {3, 0}, {4, 0}, {5, 0}, {6, 0}, {7, 0}, {8, 0}, {9, 
    0}, {10, 0}, {11, 0}, {12, 0}, {13, 0}, {14, 0}, {15, 0}, {16, 0}, {17, 
    0}, {18, 0}, {19, 0}, {20, 0}, {21, 0}, {22, 0}, {23, 0}, {24, 0}, {25, 
    0}, {26, 0}, {27, 0}, {28, 0}, {29, 0}, {30, 0}, {31, 0}, {32, 0}, {33, 
    0}, {34, 0}, {35, 0}, {36, 0}, {37, 0}, {38, 0}, {39, 0}, {40, 0}, {41, 
    0}, {42, 0}, {43, 0}, {44, 0}, {45, 0}, {46, 0}, {47, 0}, {48, 0}, {49, 
    0}, {50, 0}}
]
\dataplot[showpoints=true]{\moreverts}
\psline{<->}(0,35)(0,-10)
\psline{*-*}(0,0)(0,31)
\uput[d](24.5,0){$e_{25}$}
\uput[d](50,0){$v_{50}$}
\uput[l](0,31){$31$}
\end{pspicture}
\end{center}

The matrices and Gaussian type plots we have encountered are reminiscent of 
those that appear in the study of discrete differential operators.
We emphasize here that this phenomena results from the inner product
or metric, in particular its inverse.

From our computation of $\bigstar$ one can easily compute $\bigstar^2$, 
and it is clear that this operator approximates a delta-type function.

\bigskip

{\sc Scott O. Wilson, Department of Mathematics, Stony Brook University,
Stony Brook, NY 11796.}

email: {\tt wilson@math.sunysb.edu}


\begin{thebibliography}{99}
\bibitem{DA} Adams, D. ``A Doubled Discretisation of Chern-Simons Theory,''
arxiv.org:hep-th/9704150. 
\bibitem{BS} de Beauc\'e, V, and Sen, Samik ``Chiral Dirac Fermions on the 
Lattice using Geometric Discretisation,'' arxiv.org:hep-th/0305125.
\bibitem{BS2} de Beauc\'e, V, and Sen, S ``Discretizing Geometry and 
Preserving Topology I: A Discrete Calculus,'' arxiv.org:hep-th/0403206. 
\bibitem{BR} Birmingham, D. and Rakowski, M. ``A Star Product in Lattice
Gauge Theory,'' Phys. Lett. B 299 (1993), no.3-4, 299-304.
\bibitem{CS} Cheeger, J. and Simons, J ``Differential Characters and Geometric 
Invariants,'' Lecture Notes in Mathematics, no. 1167, 50-81, Springer.
\bibitem{Do} Dodziuk, J. ``Finite-Difference Approach to the Hodge Theory
of Harmonic Forms,'' Amer. J. of Math. 98, No. 1, 79-104.
\bibitem{DP} Dodziuk J. and Patodi V. K. ``Riemannian Structures and 
Triangulations of Manifolds,'' Journal of Indian Math. Soc. 40 (1976) 1-52.
\bibitem{Dper} Dodziuk, J. personal communication.
\bibitem{JD} Dupont, J. ``Curvature and Characteristic Classes,'' Lecture Notes
in Mathematics, vol. 640, Springer-Verlag 1978.
\bibitem{EC} Eckmann, B. ``Harmonische Funktionnen und Randvertanfgaben
in einem Komplex,'' Commentarii Math. Helvetici, 17 (1944-45), 240-245.
\bibitem{FK} Farkas, H. and Kra, I. ``Riemann Surfaces,'' Springer-Verlag,
1991. 
\bibitem{GK} Gross, P. and Kotiuga, P. R. ``Electromagnetic Theory and 
Computation: a Topological Approach,'' Cambridge University Press, Cambridge, 
2004. x+278 pp.
\bibitem{GY} Gu, X. and Yau, S.T. ``Computing Conformal Structures of
Surfaces,'' Comm. Inf. Sys. 2 (Dec. 2002) no.2, 121-146.
\bibitem{JJ} Jin, J. ``The Finite Element Method in Electrodynamics''
(Wiley, NY, 1993).
\bibitem{RK} Kotiuga, R. ``Hodge Decompositions and Computational 
Electrodynamics,'' Ph.D. thesis (McGill U., Montreal, Canada) 1984.
\bibitem{MK} Kervaire, M. ``Extension d'un theorem de G de Rham et 
expression de l'invariant de Hopf par une integrale,'' C. R. Acad. Sci. Paris,
237 (1953) 1486-1488.
\bibitem{YM} Manin, Y. ``The Partition Function of the Polyakov String can
be Expressed in Terms of Theta-Functions,'' Phys. Lett. B 172 (1986), no. 2,
184-185.
\bibitem{MS} Costa-Santos, R. and McCoy, B.M. ``Finite Size Corrections for 
the Ising Model on Higher Genus Triangular Lattices,'' J. Statist. Phys.
112 (2003), no.5-6, 889-920. 
\bibitem{CM} Mercat, C. ``Discrete Riemann Surfaces and the Ising Model,''
Comm. Math. Phys. 218 (2001), no. 1, 177-216.
\bibitem{CM2} Mercat, C. ``Discrete period Matrices and Related Topics,''
arxiv.org math-ph/0111043, June 2002.
\bibitem{CM3} Mercat, C. ``Discrete Polynomials and Discrete Holomorphic 
Approximation,'' arxiv.org math-ph/0206041.
\bibitem{RS} Ranicki, A. and Sullivan, D. ``A Semi-local Combinatorial
Formula for the Signature of a $4k$-manifold,'' J. Diff. Geometry, Vol, II 
(1976), p23-29.
\bibitem{SSSA} Sen, Sen, Sexton and Adams ``Geometric Discretisation Scheme 
Applied to the Abelian Chern-Simons Theory,'' Phys. Rev. E (3) 61 (2000), no.
3, 3174--3185.
\bibitem{SL} Smits, L. `` Combinatorial Approximation to the Divergence
of 1-forms on Surfaces,'' Israel J. of Math., vol. 75 (1991) 257-71.
\bibitem{Sp} Spivak, M. A. ``Comprehensive Introduction to Differential
Geometry,'' vol IV, Publish or Perish Inc., Boston, MA, 1975.
\bibitem{GS} Springer, G. ``Introduction to Riemann Surfaces,'' 
Addison-Wesley Publ. Company, Reading, MA, 1957. 
\bibitem{DPS} Sullivan, D. ``Infinitesimal Computations in Topology,'' 
IHES vol. 47 (1977) 269-331.
\bibitem{DS} Sullivan, D. ``Local Constructions of Infinity Structures,''
preprint.
\bibitem{TKB} Tarhassaari, T., Kettunen, L., and Bossavit, A. ``Some
Realizations of a Discrete Hodge Operator: A Reinterpretation of Finite Element
Techniques,'' IEEE Trans. Magn. vol. 35, no. 3, May 1999.
\bibitem{TC} Teixeira, F. L. and Chew, W. C. ``Lattice Electromagnetic
Theory from a Topological Viewpoint,'' J. of Math. Phys. 40 (1999) 169-187. 
\bibitem{TT} Tradler, T. and Mahmoud Zeinalian. ``Poincare Duality at the 
Chain Level, and a BV Structure on the Homology of the Free Loops Space of a 
Simply Connected Poincare Duality Space,'' arxiv math.AT/0309455.
\bibitem{Wh} Whitney, H. ``Geometric Integration Theory,'' 
Princeton Univ. Press, Princeton, NJ, 1957.
\bibitem{Wh2} Whitney, H. ``On Products in a Complex,'' Annals of Math. (2)
39 (1938), no. 2, 397-432.
\end{thebibliography}
\end{document}